\documentclass[11pt]{article}

\usepackage[ansinew]{inputenc}
\usepackage{amsmath}
\usepackage{amsfonts}
\usepackage{amssymb}
\usepackage{amsthm}
\usepackage{graphicx}
\usepackage{color}
\usepackage{multirow}
\usepackage{enumerate}
\usepackage[all]{xy}
\usepackage{url}
\usepackage{varioref}

\usepackage{xypic}

\newcommand{\N}{\mathbb{N}}
\newcommand{\Z}{\mathbb{Z}}
\newcommand{\Q}{\mathbb{Q}}
\newcommand{\C}{\mathbb{C}}
\newcommand{\ra}{\rightarrow}

\newtheorem{theorem}{Theorem}[section]
\newtheorem{proposition}[theorem]{Proposition}
\newtheorem{lemma}[theorem]{Lemma}
\newtheorem{corollary}[theorem]{Corollary}

\newtheorem{remark}[theorem]{Remark}
\newtheorem*{example}{Example}

\newcounter{HPcont}[section]
\renewcommand{\theHPcont}{(P
)}

\numberwithin{equation}{section}

\title{Quotients of the Dwork Pencil}

\author{Gilberto Bini\footnote{Universit\`{a} degli Studi di Milano - Dipartimento di Matematica ``F. Enriques'' - Via C. Saldini, 50 - 20133 Milano (Italy). {\sc E-mail}: {\tt gilberto.bini@unimi.it}}$\,$ and Alice Garbagnati\footnote{Universit\`{a} degli Studi di Milano - Dipartimento di Matematica ``F. Enriques'' - Via C. Saldini, 50 - 20133 Milano (Italy). {\sc E-mail}: {\tt alice.garbagnati@unimi.it}}}
\date{\today}

\begin{document}

\maketitle

\begin{abstract}
In this paper we investigate the geometry of the Dwork pencil in any dimension. More specifically, we study the automorphism group $G$ of the generic fiber of the pencil over the complex projective line, and the quotients of it by various subgroups of $G$. In particular, we compute the Hodge numbers of these quotients via orbifold cohomology.
\end{abstract}
\tableofcontents

\section{Introduction}\label{sec1}

The investigation of Calabi--Yau varieties is motivated by their relevance in the classification of varieties with trivial canonical bundle (as shown by the Bogomolov Decomposition Theorem) as well as in mirror symmetry. In dimension one, Calabi--Yau varieties are elliptic curves and everything is known on the moduli space and the automorphism groups of these curves. In dimension two,  Calabi--Yau varieties are K3 surfaces and several results on the automorphism groups are known, but a complete classification of the families of K3 surfaces having a given automorphism group is far from being completely understood. In higher dimension, very few results on the automorphisms are known.

If $G$ is a finite group of automorphisms of a Calabi--Yau variety $V$, and the quotient $V/G$ has a desingularization which is a Calabi--Yau variety, then every element in $G$ acts trivially on the period of $V$. Vice versa, if every element in a finite group $G\subset Aut(V)$ preserves the period of $V$, under certain assumptions the quotient $V/G$ has a desingularization which is a Calabi--Yau variety. In particular, if the dimension of $V$ is $dim(V)\leq 3$, then $V/G$ always has a desingularization which is a Calabi--Yau variety (see, for instance, \cite{bkr}). This allows one to obtain a family of Calabi--Yau varieties from another one (e.g., by taking the quotient) and in certain cases the two families are mirror: see, for instance, \cite{Bat}. This is indeed one of the most classical construction of two mirror families of Calabi--Yau 3-folds. Let $\mathcal{F}$ be the family of all quintics in $\mathbb{P}^4$, and let $X_\lambda^5$ be a special subfamily, which admits a certain group, $H_4\simeq(\Z/5\Z)^3$, as a group of automorphisms preserving the period. The one-dimensional family of Calabi--Yau 3-folds $Y_\lambda^5$ obtained by desingularizing $X_\lambda^5/H_4$ is the mirror of $X_\lambda^5$. The family $X_\lambda^5$ is known as "the Dwork pencil": see, for instance, \cite{can}.

Similarly, one can consider the one-dimensional subfamily $X_\lambda^{n+1}$ of the family of Calabi--Yau $(n-1)$-folds which are hypersurfaces of degree $n+1$ in $\mathbb{P}^n$ such that $X_\lambda^{n+1}$ admits a certain group $H_n\simeq (\Z/(n+1)\Z)^{n-1}$ as a group of automorphisms preserving the period. It turns out that $X_\lambda^{n+1}$ is  $V(\sum_{i=1}^{n+1}x_i^{n+1}-(n+1)\lambda\prod_{i=1}^{n+1}x_i)$ and it is immediate to show that the group of automorphisms of $X_\lambda^{n+1}$ which preserves the period is greater than $H_n$ and is in fact $\mathfrak{A}_{n+1}\rtimes H_n$. If the dimension of $X_\lambda^{n+1}$ is either $1$ or $2$, this phenomenon is better known and more general; indeed, every elliptic curve admits $\mathfrak{A}_{3}\rtimes H_2\simeq (\Z/3\Z)^2$ as a group of automorphisms preserving the period. Moreover, every K3 surface with $H_3:=(\Z/4\Z)^2$ as a group of automorphisms preserving the period has also the group $\mathfrak{A}_4\rtimes H_3$ as a group of automorphisms preserving the period.

The aim of this paper is to describe the automorphisms and the quotients of the family $X_\lambda^{n+1}$ in every dimension. In Section \ref{sec Dwork pencil definitions} we introduce the Dwork pencil, i.e., the Calabi--Yau $(n-1)$-folds $X_\lambda^{n+1}$. Section \ref{sec:autom group X} is the fundamental part of our paper, where we give the automorphism group of the family $X_\lambda^{n+1}$ (and of some of its special members), and we describe several quotients. In particular, we observe that the automorphism group of the abstract variety and of the embedded variety $X_\lambda^{n+1}$ coincides if and only if $n>3$. The automorphism group of the embedded $X_\lambda^{n+1}$ is $\mathfrak{S}_{n+1}\rtimes H_{n}$ and its subgroup $\mathfrak{A}_{n+1}\rtimes H_{n}$ preserves the period.

Given the automorphism group, a very natural problem is to describe the quotient varieties: sections \ref{subsec: quotient by A}, \ref{subsec: Y}, \ref{subsec: quotients by G} and \ref{subsec: quotient by tau} are devoted to the description of $X_\lambda^{n+1}/\mathfrak{S}_{n+1}$ and $X_\lambda^{n+1}/\mathfrak{A}_{n+1}$, of $X_\lambda^{n+1}/H_n$, of $X_\lambda^{n+1}/\mathfrak{S}_{n+1}\rtimes H_{n}$ and $X_\lambda^{n+1}/\mathfrak{A}_{n+1}\rtimes H_{n}$ and $X_\lambda^{n+1}/\tau$, where $\tau$ is an involution which does not preserve the period. In Section \ref{subsec: quotient by A}, we write the equations of $X_\lambda^{n+1}/\mathfrak{S}_{n+1}$ and $X_\lambda^{n+1}/\mathfrak{A}_{n+1}$ as singular subvarieties of weighted projective space. The orbifold $X_\lambda^{n+1}/\mathfrak{A}_{n+1}$ is a Calabi--Yau orbifold and it admits a desingularization which is a smooth Calabi--Yau if and only if the dimension is less than or equal to 3. In Section \ref{subsec: Y},  we describe the very classical quotient $X_\lambda^{n+1}/H_n$ which is a singular orbifold which has a Calabi--Yau variety, $Y_\lambda^{n+1}$, as resolution. We observe that if $n>4$, the two quotients $X_\lambda^{n+1}/\mathfrak{A}_{n+1}$ and $X_\lambda^{n+1}/H_{n}$ are examples of two opposite situations: in both cases they are quotient of a Calabi--Yau variety by a group of automorphisms preserving the period, but in the first case the quotient can not be desingularized to a Calabi--Yau variety, and in the second case it admits a Calabi--Yau variety as a desingularization.

Sections \ref{sec: K3 X} and \ref{sec: Quotients CY 3 fold}  are devoted to the particular cases of dimension $2$ and $3$, respectively. In Section \ref{sec: K3 X} we describe the N\'eron--Severi group of $X_\lambda^4$ and of certain special members. We use this description in order to find two lattices, called $\Omega_{\mathfrak{A}_4}$, $\Omega_{\mathfrak{S}_4}$, which are related to the presence of the groups $\mathfrak{A}_4$ and $\mathfrak{S}_4$ as group of automorphisms which preserve the period on a K3 surface. The quotient surfaces $X_\lambda^4/\mathfrak{A_4}$, $X_\lambda^4/H_3$, $X_\lambda^4/(\mathfrak{A_4}\rtimes H_3)$ are described in Proposition \ref{prop: K3 quotient A4}, Section \ref{subsub K3 Y} and Section \ref{subsub K3 quotient A4}, respectively.

In Section \ref{sec: Quotients CY 3 fold}, we consider the quotient of $X_\lambda^5$ by some subgroups of the automorphism group $\mathfrak{S}_5\rtimes H_4$ . First, we give some preliminary and general results on the quotient of a Calabai--Yau 3-fold by automorphisms preserving the period and, after that, we compute the Hodge numbers of a Calabi--Yau desingularization of $X_\lambda^5/J$, where $J\subset\mathfrak{A}_5\rtimes H_4$. The results of this Section are of a different type with respect to the ones of Section \ref{sec:autom group X}. Indeed, in Section \ref{sec:autom group X}, we give a geometrical description of the quotient; here we compute the Hodge numbers of the quotient via orbifold cohomology and without providing a geometrical construction.

\section{The Dwork Pencil $X_{\lambda}^{n+1}$}\label{sec2}\label{sec Dwork pencil definitions}

In this section, we recall the definition of the Dwork pencil. Let $n \geq 2$ be a positive integer. Denote by $X_{\lambda}^{n+1}$ the zero locus $Z(F_{\lambda}^{n+1})$ in complex projective space ${\mathbb P}^{n}$, where
$$
F_{\lambda}^{n+1}:=\sum_{i=1}^{n+1}x_i^{n+1} - (n+1)\lambda \prod_{j=1}^{n+1}x_j \in {\mathbb C}[x_1, \ldots, x_{n+1}],
$$
and $\lambda$ is a complex number. To begin with, we investigate singularities of $X_{\lambda}^{n+1}$. From now on, $\xi_s$ is a primitive $s$-th root of unity.

\begin{proposition}
i) The algebraic variety $X_{\lambda}^{n+1}$ is smooth for any $\lambda$ such that $\lambda^{n+1} \neq 1$.\\
ii) For $\lambda$ such that $\lambda^{n+1}=1$, the variety
$X_{\lambda}^{n+1}$ has $(n+1)^{n-1}$ ordinary double  points. For $\lambda=\xi_{n+1}^r$ the
singular points are given by
\begin{equation}
\label{first}
\left\{(\xi_{n+1}^{i_1}: \ldots: \xi_{n+1}^{i_n}: 1): \sum_j i_j \equiv -r \,
\, \hbox{mod} \, \, (n+1) \right \}.
\end{equation}
\\
\end{proposition}

\proof i) The Jacobian ideal is given by
$$
\left( (n+1)x_1^n - (n+1) \lambda \prod_{j \neq 1} x_j, \ldots,
(n+1) x_{n+1}^{n} - (n+1)\lambda\prod_{j=1}^nx_j\right).
$$

After multiplying the equation $x_k^n - \lambda \prod_{j \neq k}x_j=0$ by $x_k$ for $k=1, 2, \ldots, n+1$, we obtain the following identity:
\begin{equation}
\label{prodotto} \prod_{k} x_k^{n+1} = \lambda^{n+1} \prod_k
x_k^{n+1}.
\end{equation}

If one of the coordinates $x_j$ is zero, and the first order partial derivatives are zero, all coordinates are zero. Thus, \eqref{prodotto} yields $\lambda^{n+1}=1$.

ii) As noticed before, $x_{n+1}$ can not be zero. Set $x_{n+1}=1$,
so the equation $x_{n+1}^n - \lambda x_1 \cdots x_n=0$ yields
$x_1\cdots x_n=1/\lambda$ for $\lambda \neq 0$.  By plugging in
the other equations (given by the first order partial derivatives),
we get $x_k^{n+1}=1$ for any $k=1, \cdots, n$. In other
words, $x_k = \xi_{n+1}^{i_k}$ for some $i_k \in {\mathbb Z}/ (n+1)
{\mathbb Z}$. Since $x_1 \ldots x_n=\xi_{n+1}^{-r}$ has to be satisfied, \eqref{first} is proved.  Notice that there exist transformations which map the nodes
one onto the other. In fact, if $P=(\xi_{n+1}^{i_1}, \xi_{n+1}^{i_2}, \ldots ,
1)$ and $Q= (\xi_{n+1}^{j_1}, \ldots , 1)$ are two nodes, the map $y_k=
\xi_{n+1}^{j_k - i_k}x_k$ for $k=1, \ldots, n$ and $y_{n+1}=x_{n+1}$
will do. Thus, it suffices to study the type of the singularity of
$(1, 1, 1, \ldots , 1)$. This is done in \cite{schoen}.
\endproof

\begin{remark}For $\lambda^{n+1}\neq 1$ the varieties
$X^{n+1}_{\lambda}$ are smooth hypersurfaces of degree $n+1$ in
$\mathbb{P}^n$, hence they are Calabi-Yau manifolds.\end{remark}

\subsection{The Hodge Numbers of $X_{\lambda}^{n+1}$}\label{sec3}\label{subsec: Hodge numbers X}

Let $S$ be a smooth hypersurface in $\mathbb{P}^{m+1}$ of degree $d$. By the Lefschetz Hyperplane Theorem, we deduce that $h^{i}(S)=1$ for even $i$ such that $0\leq i\leq 2m$ and $i\neq m$, as well as $h^i(S)=0$ for odd $i$ and $i\neq m$. More precisely, by the Hodge decomposition of $H^i(S,\C)$ for $0\leq i\leq 2m$ and $i\neq m$, the group $H^{i/2,i/2}(S,\C)$ is isomorphic to $\C$ for even $i \neq m$. Thus, all the other Hodge groups of weight $i\neq m$ are trivial.

In order to compute the Hodge diamond of $S$, it suffices to find the Hodge decomposition of the middle cohomology. This can be done by the following well known formula, due to Griffiths (see, for instance, \cite{V1}):

\begin{equation}\label{formula: Griff hodge nubers}H^{p,m-p}_0(S)=\left(\C[x_1,\ldots x_{m+2}]/J_F\right)_{d(p+1)-(m+2)},\end{equation}
where $S=V(F)$ and $H^{p,m-p}_0(S)$ is the primitive cohomology. We recall that $h^{p,m-p}_0(S)=h^{p,m-p}(S)$ for every $p$ and any odd $m$. If $m$ is even, then $h^{p,m-p}_0(S)=h^{p,m-p}(S)$ for $p\neq m/2$ and $h^{p,m-p}_0(S)=h^{p,m-p}(S)-1$ for $p=m/2$. If we are interested in the Betti numbers (and not in the Hodge ones), the computation is certainly easier: indeed, it is well known that the Euler characteristic of $S$ is given by
\begin{equation}\label{Euler char}e(S)=\sum_{k=0}^m(-1)^kd^{k+1}\binom{m+2}{m-k}.\end{equation}

Therefore, we have
\begin{equation}
\label{hm}h^m(S)=(m+1)-e(S)\, \hbox{ if } \, m \, \hbox{is odd},
\end{equation}
\begin{equation}
h^m(S)=e(S)-m \, \hbox{ if }\,  m \, \hbox{is even}.
\end{equation}

\begin{example}{\rm
For $n=3$, the variety $X_\lambda^4$ is a K3 surface and its Hodge numbers are uniquely determined by this property:
$$h^{0,0}=h^{2,0}=h^{2,2}=1, h^{1,0}=0, h^{1,1}=20.$$

Also, the Hodge numbers can be computed explicitily because $X_{\lambda}^4$ is a Calabi-Yau variety, so one has to compute $h^{1,1}$. By \eqref{Euler char}, the Euler characteristic is $e(X_\lambda^4)=24$; hence, $h^{2}(X_\lambda^4)=22$ by \eqref{hm}, which yields $h^{1,1}(X_\lambda^4)=20$.

For $n=4$, it suffices to compute the Hodge numbers of the middle cohomology for the Hodge diamond of $X_\lambda^5$. The only numbers one has to determine are $h^{2,1}$ and $h^{1,2}$, which coincide. Again, it suffices to apply Formulae \ref{Euler char} and \ref{hm} in order to conclude that $h^{2,1}=101$.

For $n=5$, it suffices to calculate the Hodge numbers $h^{3,1}$ and $h^{2,2}$. By \eqref{Euler char} and \eqref{hm}, we obtain $e(X_\lambda^6)=2610$ and $h^4(X_\lambda^6)=2606$, so one has $$h^{3,1}(X_\lambda^6)+h^{2,2}(X_\lambda^6)+h^{1,3}(X_\lambda^6)=2604,$$
which is not enough to determine the whole Hodge diamond. For this purpose, we apply Formula \ref{formula: Griff hodge nubers}. Since all the hypersurfaces of a given dimension and degree have the same Hodge numbers, we can
chose $F$ to be the polynomial $\sum_{i=1}^6x_i^6$; hence the Jacobian ideal is generated by $x_i^5$ for $i=1,\ldots, 6$. This implies that
 $$h^{1,3}(X_\lambda^6)=\dim\left(\left(\C[x_1,\ldots x_{6}]/(x_1^5,\ldots x_6^5)\right)_{6}\right).$$

Now, it is clear that the space $\left(\C[x_1,\ldots x_{6}]/(x_1^5,\ldots x_6^5)\right)_{6}$ is generated by the monomials $\prod_{i=1}^6 x_i^{\alpha_i}$ such that $\sum_{i=1}^6\alpha_i=6$ and $\alpha_i<5$ for each $i=1,\ldots 6$. The dimension of this space is $\binom{6+5}{6}-6^2=426$. Hence, the non trivial Hodge numbers of $X_\lambda^6$ are $$h^{0,0}=h^{1,1}=h^{4,0}=h^{0,4}=1,\,\,  h^{3,1}=h^{1,3}=426,\, \,  h^{2,2}=1752.$$
}
\end{example}

\begin{remark}{\it More generally, let $V$ be a smooth hypersurface of degree $m+2$ in $\mathbb{P}^{m+1}$ (i.e., a Calabi--Yau hypersurface). By the same arguments as before, i.e., by Formula \ref{formula: Griff hodge nubers} and by choosing $F$ to be the Fermat polynomial, the following holds:  $$h^{1,m-2}=\binom{2m+3}{m+2}-(m+2)^2,$$ which is the dimension of the family of Calabi--Yau manifolds that are hypersurfaces of degree $m+2$ in $\mathbb{P}^{m+1}$}.
 \end{remark}

\section{The Automorphism Group of $X_{\lambda}^{n+1}$}\label{sec7}\label{sec:autom group X}

Let us recall that a Calabi-Yau manifold $V$ of complex dimension
$n-1$ admits a nowhere vanishing holomorphic $(n-1)$-form, i.e.,
$H^{n-1,0}(V)=\langle\omega_V\rangle$, which is called the {\em period} of $V$.

Let us consider the following action of $({\Z/(n+1)\Z})^{n-1}$ on
$\mathbb{P}^n$. Denote by $\xi_{n+1}$ a primitive $(n+1)$-th root of unity
and define automorphisms of ${\mathbb P}^n$ as follows:
$$
h_{(a_1, \ldots, a_{n+1})}(x_1: \ldots : x_{n+1}):=(\xi_{n+1}^{a_1}x_1:
\ldots : \xi_{n+1}^{a_{n+1}}x_{n+1}).
$$
Consider the finite group
\begin{equation}\label{formula: group H}
H_n= \langle h_1:=h_{(1,-1,\ldots, 0)}, \ldots, h_{n-1}:=h_{(1,
0, \ldots, -1, 0)} \rangle \cong \left( {\mathbb Z}/(n+1){\mathbb
Z}\right)^{n-1}.
\end{equation}
We notice that $H_n$ is the subgroup of the group generated by the automorphisms $h_{(a_1,\ldots a_{n+1})}$ obtained requiring $a_1=1$ and $\sum_i a_i\equiv 0\mod (n+1)$.\\

Now, fix a generic $\lambda \in {\mathbb C}$.
\begin{proposition}\label{prop: Z/(n+1)Z^(n-1) group of autom}
The family $X_{\lambda}^{n+1}$ is the family of the hypersurfaces
of degree $n+1$ in $\mathbb{P}^n$ which admit the group $H_n$ as
a group of automorphisms. The automorphisms of $H_n$ fix the period of $X_{\lambda}^
{n+1}$.
\end{proposition}
\proof If a hypersurface of degree $n+1$ in $\mathbb{P}^n$ admits
$H_n$ as group of automorphisms, then it is a member of the family
\begin{equation}\label{formula: equation
Pn}\sum_{i_1}^{n+1}\alpha_ix_i^n+\beta
\prod_{i=1}^{n+1}x_i=0.\end{equation} The equation depends on $n+2$
parameters. The automorphisms of $\mathbb{P}^{n}$ commuting with those of $H_n$ belong to the group $GL(1)^{n}/GL(1)$. Hence, this
family has $(n+2-1)-(n+1-1)=1$ moduli. Up to a projectivity, we can
assume that the hypersurfaces in \eqref{formula: equation Pn} have
the equation $\sum_{i=1}^{n+1}x_i^{n+1}-\lambda
\prod_{i=1}^{n+1}x_i=0$ - it suffices to consider the projectivity
$diag(1/\alpha_{i}^{1/{n+1}})$. This is the equation of the
family $X^{n+1}_{\lambda}$.

The holomorphic $(n-1)$-form of a hypersurface of
$\mathbb{P}^{n}$ with equation $F_{n+1}=0$ is given by
\begin{eqnarray}\label{formula: holomorphic form}Res\left(\frac{\sum_{i=1}^{n}(-1)^ix_idx_1\wedge
dx_2\wedge\ldots\wedge\hat{dx_i}\wedge\ldots\wedge
dx_n)}{F_{n+1}}\right).\end{eqnarray} Since the equation defining
$X_{\lambda}^{n+1}$ is invariant under $H_n$, the
automorphisms in $H_n$ fix the period of $X^{n+1}_{\lambda}$.\endproof

By Proposition \ref{prop: Z/(n+1)Z^(n-1) group of autom}, the hypersurfaces $X_{\lambda}^{n+1}$
are characterized by the fact that they admit $H_n$ as a group of
automorphisms fixing the period. The next proposition shows that they
in fact admit a larger group of automorphisms, i.e.,
the hypersurfaces of degree $n+1$ in $\mathbb{P}^n$ having
$H_n\simeq (\Z/(n+1)\Z)^{n-1}$ as a group of automorphisms also have the group $\mathfrak{S}_{n+1}$ as a group of automorphisms, and in particular ${\mathfrak A}_{n+1}$ (the alternating group of
degree $n+1$) as a group of automorphisms fixing the period.

\begin{proposition}\label{prop: An and Sn groups of autom} If $n\geq 4$ and $\lambda$ is generic, then $\mathfrak{S}_{n+1}\rtimes H_n$ is the automorphism group of the abstract variety $X_\lambda^{n+1}$ and $\mathfrak{A}_{n+1}\rtimes H_n$ is the subgroup fixing the period of $X_\lambda^{n+1}$.\\
If $n=2,3$ and $\lambda$ is generic, then $\mathfrak{S}_{n+1}\rtimes H_n\subset Aut(X_\lambda^n)$ is the subgroup of the automorphism group fixing the polarization which embeds $X_\lambda^{n+1}$ in $\mathbb{P}^n$ as the zero locus $Z(F_\lambda^{n+1})$. Moreover, its subgroup $\mathfrak{A}_{n+1}\rtimes H_n$ acts on $X_\lambda^{n+1}$ fixing the period. \end{proposition}
\proof Let $\sigma\in {\mathfrak S}_{n+1}$. Then $\sigma$ defines the projectivity
$$\sigma_{\mathbb{P}^{n+1}}:(x_1:\ldots :x_{n+1})\ra (x_{\sigma(1)}:\ldots:
x_{\sigma(n+1)}).$$ The equation of $X_{\lambda}^{n+1}$ is
invariant under the action of ${\mathfrak S}_{n+1}$ and the form
defined in \eqref{formula: holomorphic form} is invariant under ${\mathfrak A}_{n+1}$. \\
Let $G$ be the group of automorphisms of $X_\lambda^{n+1}$, which preserve the embedding - given by the polarization - of $X_{\lambda}^{n+1}$ in $\mathbb{P}^{n}$. If we restrict to a hyperplane and apply known results in \cite{Sh}, the group $G$ is isomorphic to the semidirect product of the symmetric group ${\mathfrak S}_{n+1}$ and the group $H_n \cong (\Z/(n+1)\Z)^{n-1}$, where $H_n$ is normal in the semidirect product.  Thus, the following exact sequence holds:
$$
1 \rightarrow \mathfrak{A}_{n+1}\rtimes H_n \rightarrow \mathfrak{S}_{n+1}\rtimes H_n \rightarrow <\tau> \rightarrow 1,
$$
where $\tau$ is the transposition of the coordinates that swaps $x_1$ and $x_2$. Clearly, $\tau$ does not preserve the period of $X_{\lambda}^{n+1}$. The subgroup $ \mathfrak{A}_{n+1}\rtimes H_n$ is the group of automorphisms which preserve the period of $X_\lambda^{n+1}$. \\
If $n\geq 4$, $H^{1,1}(X_\lambda^n)\simeq \C$ and $Pic(X_\lambda^n)\simeq \langle H\rangle $ for a certain ample line bundle $H$. Since an automorphism of the abstract variety $X_{\lambda}^{n+1}$  maps an ample line bundle - with a given self intersection - to an ample line bundle with the same self intersection, every automorphism of $X_\lambda^{n+1}$ sends $H$ to $H$, and thus preserves the polarization associated to the embedding of $X_\lambda^{n+1}$ in $\mathbb{P}^n$.\endproof

\begin{remark} Clearly, if $n=2$ the automorphism group is greater then $\mathfrak{S}_3\rtimes H_3$, since $X_\lambda^3$ is an elliptic curve defined over $\C$ and thus its automorphism group is infinite. Similarly the automorphism group of $X_\lambda^4$ is greater then $\mathfrak{S}_4\rtimes H_4$, and in fact it is infinite. This follows from  \cite[Theorem]{kondo}. Let $G$ be one of the following groups: $\{1\}$, $\Z/2\Z$, $(\Z/2\Z)^2$, $\mathfrak{S}_3\times \Z/2\Z$ and let $S$ be a K3 surface with a finite automorphism group $Aut(X)$. Then $Aut(X)$ is either one of the groups $G$ or a cyclic extension of one of the groups $G$. Since $\mathfrak{S}_4\rtimes H_4$ is neither a group $G$ nor a cyclic extension of $G$, the group $Aut(X_\lambda^4)$ is infinite.\end{remark}

\begin{remark}\label{rem: autom X_0} For specific values of $\lambda$ the automorphism group of the abstract variety $X_\lambda^{n+1}$ might be greater then $\mathfrak{S}_n\rtimes H_n$, also for $n\geq 4$. For example, if $\lambda=0$, the group $(\Z/(n+1)\Z)^n(\supset H_n)$ generated by the $h_{(a_1,\ldots a_{n+1})}$ such that $a_1=1$, is a group of automorphisms of $X_0^{n+1}$. In this case, also the automorphism group preserving the period might be larger. For example, if $\lambda=0$ and $n\equiv 1\mod 2$, then $\mathfrak{S}_{n+1}$ is a subgroup of the automorphism group fixing the period: indeed, the composition of $\tau$ and $h_{(0,0,0,\ldots -1)}$ is an automorphism of $X_0^{n+1}$ preserving the form \eqref{formula: holomorphic form} \end{remark}

We observe that the family of hypersurfaces of degree $n+1$ in
$\mathbb{P}^n$ admitting ${\mathfrak S}_{n+1}$ as a group of automorphisms and
the one of hypersurfaces of degree $n+1$ in $\mathbb{P}^n$
admitting $H_n\simeq (\Z/(n+1)\Z)^{n-1}$ as a group of automorphisms
are not the same, in particular the former specializes to the latter. Indeed, the family of hypersurfaces of degree
$n+1$ in $\mathbb{P}^n$ admitting ${\mathfrak S}_{n+1}$ as a group of automorphisms is constructed as zeros of a homogeneous symmetric polynomial $f(x_1:\ldots : x_{n+1})$ of degree $n+1$. In particular, $f$ can be written as a linear combination of $p(n+1)$ linearly independent terms, where $p(n+1)$ is the partition function of $n+1$. It is well known that $p(n+1)\geq (n+1)$, hence the polynomial $f$ depends at least on $(n+1)$ parameters. The group of projectivities commuting with ${\mathfrak S}_{n+1}$ is made up of the matrices of the form $A=[a_{i,j}]$, where $a_{i,j}=1$ if $i\neq j$, and $a_{i,i}=c$ with $c\in \C$. This implies that the family of hypersurfaces of degree
$n+1$ in $\mathbb{P}^n$ admitting ${\mathfrak S}_{n+1}$ as a group of automorphisms depends at least on $(n+1)-1-1=n-1$ parameters. For each $n>2$, this family is bigger then the family of hypersurfaces of degree $n+1$ in $\mathbb{P}^{n}$ admitting $H_n$ as a group of automorphisms. Therefore, they cannot coincide.

Notably, the quotient of a Calabi--Yau variety by a finite group which preserves the period may admit a desingularization which is a smooth Calabi--Yau variety. In fact, the desingularization certainly exists if the dimension of the Calabi--Yau is less than 4. On the contrary, the quotient of a Calabi--Yau variety, which does not preserve the period, can not have a desingularization which is a Calabi--Yau variety, because it does not have any top-degree holomorphic form.

\subsection{The Quotient by ${\mathfrak S}_{n+1}$ and ${\mathfrak A}_{n+1}$}\label{sec8}\label{subsec: quotient by A}

The group ${\mathfrak S}_{n+1}$ acts on $X_{\lambda}^{n+1}$ via permutation of the coordinates and its subgroup ${\mathfrak A}_{n+1}$ acts
on $X_{\lambda}^{n+1}$ by preserving the period: see Proposition \ref{prop: An and Sn groups of autom}. Here we consider the quotient of
$X_{\lambda}^{n+1}$ by these groups. Note that $\lambda$ is generic.

Recall that the quotient of projective space ${\mathbb P}^n$ by ${\mathfrak S}_{n+1}$ is isomorphic to $W{\mathbb P}^n(1,2, \ldots, n+1)$. The coordinates in $W{\mathbb P}^n(1,2, \ldots, n+1)$ are given by weighted homogeneous variables $e_i$ of weight $i$, respectively, which are the symmetric function of degree $i$ in the variables $x_j$, $j=1,\ldots n+1$. Let $w_{n+1}$ be a variable of weight $\frac{n(n+1)}{2}$ and denote by $W{\mathbb P}^{n+1}(1, 2, \ldots, n+1, \frac{n(n+1)}{2})$ the weighted projective space with weighted homogeneous coordinates $[e_1, \ldots, e_{n+1}, w_{n+1}]$. Let
$$
\pi: W{\mathbb P}^{n+1}\left(1, 2, \ldots, n+1, \frac{n(n+1)}{2}\right) \rightarrow W{\mathbb P}^n(1,2, \ldots, n+1)
$$
be the projection onto the first $n+1$ coordinates. Let
$$
A_{n+1}:=\left\{[e_1, \ldots, e_n, e_{n+1}, w_{n+1}] \in W{\mathbb P}^{n+1}\left(1, 2, \ldots, n+1, \frac{n(n+1)}{2}\right)
\right.
$$
$$
\left.
 : w_{n+1}^2 - \theta_{n+1} w_{n+1} - \Delta_{n+1}=0\right\},
$$
where $\theta_{n+1}$ and $\Delta_{n+1}$ are the symmetric polynomials
$$
\prod_{1 \leq i < j \leq n+1} (x_i + x_j), \quad \frac{1}{4}\left( \prod_{1 \leq i < j \leq n+1} (x_i-x_j)^2 -\prod_{1 \leq i < j \leq n+1} (x_i + x_j)^2\right),
$$
respectively. Note that
the discriminant hypersurface is given by
\begin{equation}
\label{disc} \prod_{1 \leq i < j \leq n+1} (x_i - x_j)^2=det(M_{n+1})=:V^2,
\end{equation}
where
$$
M_{n+1}:=\left(
\begin{array}{ccccc}
n+1 & \psi_1 & \psi_2 & \hdots & \psi_{n} \\
\psi_1 & \psi_2 & \psi_3 & \hdots & \psi_{n+1} \\
\psi_2 & \psi_3 & \psi_4 & \hdots  &\psi_{n+2}\\
\vdots & \vdots & \vdots & \vdots & \vdots \\
\psi_{n} & \psi_{n+1} & \psi_{n+2} & \hdots & \psi_{2n}
\end{array}
\right),
$$
(see, for instance, \cite{proc}), where $\psi_j$ are the Newton power sums. As well known, $A_{n+1}$ is isomorphic to ${\mathbb P}^n/{\mathfrak A}_{n+1}$. Clearly, the projection $\pi$ restricted to $A_{n+1}$ is a degree $2$ map onto $W{\mathbb P}^n(1,2, \ldots, n+1)$.

The image of $X_{\lambda}^{n+1}$ in ${\mathbb P}^n/{\mathfrak S}_{n+1}$ is
given by $\psi_{n+1}-(n+1)\lambda e_{n+1}$. This hypersurface $S_{\lambda}^{n+1}$ is an orbifold in $W{\mathbb P}^n(1,2, \ldots, n+1)$ for any $n \geq 3$. Since the size of the symmetric group is quite large, it is more convenient to study the singularities of $S_{\lambda}^{n+1}$  in another way. First, recall that the singular locus $\Sigma_{n+1}$ of $W{\mathbb P}^n(1,2, \ldots, n+1)$ is the union of $Sing_p$ for any prime $p$, where
$$
Sing_p:=\left\{[e_1, \ldots e_{n+1}] : e_j=0 \, \, \hbox{for} \, \, p \not | j\right\}.
$$

If  $S_{\lambda}^{n+1}$ is quasi-smooth and well-formed, then the singular locus of $S_{\lambda}^{n+1}$ is the intersection of $\Sigma_{n+1}$ and $S_{\lambda}^{n+1}$ . The hypersurface $S_{\lambda}^{n+1}$  is quasi-smooth because the derivative with respect to $e_{n+1}$ does not vanish. Indeed, the power sum $\psi_{n+1}$ is given by $e_{n+1}$ and other terms. This means that the derivative of $\psi_{n+1}-(n+1)\lambda e_{n+1}$ with respect to $e_{n+1}$ is different from $0$ for generic $\lambda$. Now, we show that $S_{\lambda}^{n+1}$  is also well-formed. For any prime $p$ and any quasi-smooth weighted complete intersection $X_{d_1, \ldots, d_c}$ in $W{\mathbb P}^{n}(a_1, \ldots, a_n, a_{n+1})$ set
$$
m(p):=|\{i: p|a_i\}|, \quad k(p):=|\{i : p |d_i\}|, \quad q(p):=n-c+1 - m(p)+k(p).
$$

As proved in \cite{Di}, $X_{d_1, \ldots, d_c}$ is well-formed if and only if $q(p) \geq 2$ for all primes $p$. In our case $c=1$ and $a_i=i$ for $i=1, \ldots, n+1$.
\begin{lemma}\label{wellformed}
The quasi-smooth hypersurface $S_{\lambda}^{n+1}$  in $W{\mathbb P}^n(1,2, \ldots, n+1)$ is well-formed if and only if $n \geq 5$.
\end{lemma}

\proof Set $n+1=sp+r$, where $0 \leq r < p$. Then $q(p)=n-1-\rfloor \frac{n+1}{p} \lfloor + \varepsilon = sp+r-2-s+\varepsilon$, where $\varepsilon$ is $1$ if $p|(n+1)$ and $0$ otherwise. If $p=2$, we get $q(2)=s-1$ which is greater than or equal to $2$ if and only if $n \geq 5$. As for the other primes, an analogous computation yields the desired result.
\endproof

Take now the subvariety $T_{\lambda}^{n+1}$ in $W{\mathbb P}^{n+1}\left(1, 2, \ldots, n+1, \frac{n(n+1)}{2}\right)$ given by the equations
\begin{equation}\label{eqs}
w_{n+1}^2-w_{n+1}\theta_{n+1} - \Delta_{n+1}=0, \quad \psi_{n+1}-(n+1)\lambda e_{n+1}=0.
\end{equation}

\begin{lemma}
\label{Torbifold}
 Let $n\geq 5$. The quotient $X_{\lambda}^{n+1}/{\mathfrak A}_{n+1}$ is a Calabi-Yau orbifold in $W{\mathbb P}^n\left(1,2, \ldots, n+1, \frac{n(n+1)}{2}\right)$, which is isomorphic to $T_{\lambda}^{n+1}$. Moreover, $T_{\lambda}^{n+1}$ is a degree $2$ covering of $S_{\lambda}^{n+1}$.
\end{lemma}
\proof Clearly, $T_{\lambda}^{n+1}$ is isomorphic to the quotient $X_{\lambda}^{n+1}/{\mathfrak A}_{n+1}$and is a covering of $S_{\lambda}^{n+1}$. Moreover, $T_{\lambda}^{n+1}$ is quasi-smooth and well-formed. The former claim can be proved as follows. Take the derivative $2w_{n+1}- \theta_{n+1}$ of $w_{n+1}^2-w_{n+1}\theta_{n+1} - \Delta_{n+1}$ with respect to $w_{n+1}$. If this derivative is non zero at a point $p$, then the Jacobian of \eqref{eqs} has a two by two minor that is non zero at $p$, namely the minor corresponding to the derivatives with respect to $w_{n+1}$ and $e_{n+1}$.  Assume, now, $2w_{n+1}-\theta_{n+1}=0$. The derivative $2w_{n+1}- \theta_{n+1}$ of $w_{n+1}^2-w_{n+1}\theta_{n+1} - \Delta_{n+1}$ with respect to $e_j$ equals $-\frac 14 \frac{\partial V^2}{\partial e_j}$, where $V^2=det(M_{n+1})$. We must prove that for a point $q$ in the quotient such that $2w_{n+1}-\theta_{n+1}=0$, the matrix
$$
J:=\left(
\begin{array}{cccc}
0 & -\frac 14 \frac{\partial V^2}{\partial e_1} & \ldots & -\frac 14 \frac{\partial V^2}{\partial e_{n+1}} \\
0 & -\frac 14 \frac{\partial( \psi_{n+1}-(n+1)\lambda e_{n+1})}{\partial e_1} & \ldots &\frac{\partial( \psi_{n+1}-(n+1)\lambda e_{n+1})}{\partial e_{n+1}}
\end{array}
\right)
$$
has a two by two minor different from zero. On the contrary, suppose that all minors are zero. In particular, this implies that
$$
\frac{\partial V^2}{\partial e_n}= \mu \frac{\partial V^2}{\partial e_{n+1}},
$$
$$
\frac{\partial(( \psi_{n+1}-(n+1)\lambda e_{n+1})}{\partial e_n}= \mu \frac{\partial(( \psi_{n+1}-(n+1)\lambda e_{n+1})}{\partial e_{n+1}}.
$$

Passing to the variables $\psi_n$ and $\psi_{n+1}$, we have
$$
\frac{\partial V^2}{\partial \psi_n}= \sum_{k \geq 1}\psi_1^k F_k(\psi_2, \ldots, \psi_{n+1}).
$$

As readily checked, the polynomial $V^2$ contains the term $\psi_n^{n+1}$: its derivative with respect to $\psi_n$ is not a multiple of $\psi_1^k$ for some $k  \geq 1$. This means that the minor corresponding to the variables $e_n, e_{n+1}$ can not be zero.  This shows that $T_{\lambda}^{n+1}$ is quasi-smooth.

Analogously to Lemma \ref{wellformed}, one can prove that $T_{\lambda}^{n+1}$ is well-formed. The canonical sheaf of $T_{\lambda}^{n+1}$ is  trivial since the sum of the weights of projective space equals the sum of the degrees of the equations defining $T_{\lambda}^{n+1}$. By Proposition 4.1.3 in \cite{CK}, it suffices
to show that the singularities of $X_{\lambda}^{n+1}/{\mathfrak
A}_{n+1}$ are canonical. In fact, they are normal because they are finite quotient singularities. Furthermore, they are Gorenstein. Since the action of ${\mathfrak A}_{n+1}$ preserves the period, any element of  ${\mathfrak A}_{n+1}$  corresponds to an element in $SL(n,{\mathbb C})$. By \cite{Di}, Proposition 2, the singularities of $T_{\lambda}^{n+1}$ are those of the ambient weighted projective space. These are canonical since Gorenstein toric varieties have at worst canonical singularities. The Hodge numbers $h^{0,j}$ - for $0 < j < n-2$ - of the quotient are easily seen to be zero. Thus, the claim follows.
\endproof

For $n=3,4$ the singularities of $S_{\lambda}^{n+1}$ and $T_{\lambda}^{n+1}$ must be computed directly. First, let us focus on the case $n=3$.

\begin{proposition}\label{prop: K3 quotient A4} If $\lambda^4\neq 1$, the surface
$X_\lambda^4/\mathfrak{S}_4$ has 5 singular points of type $A_2$
and 3 singular points of type $A_1$ and the surface
$X_\lambda^4/\mathfrak{A}_4$ has 6 singular points of type $A_2$
and 4 singular points of type $A_1$.\\
If $\lambda^4= 1$, the surface $X_\lambda^4/\mathfrak{S}_4$ has 5
singular points of type $A_2$ and 6 singular points of type $A_1$
and the surface $X_\lambda^4/\mathfrak{A}_4$ has 6 singular points
of type $A_2$ and 8 singular points of type
$A_1$.\end{proposition} \proof To describe the singularities of
$X_{\lambda}^4/\mathfrak{S}_4$ (resp. of
$X_{\lambda}^4/\mathfrak{A}_4$) we consider the orbit of the
points of $X_{\lambda}^4$ with a non trivial stabilizer with
respect to
$\mathfrak{S}_4$ (resp. $\mathfrak{A}_4$).\\
Let us consider the action of $\mathfrak{A}_4$. There are 6 orbits
of points with stabilizer isomorphic to $\Z/3\Z$, the orbits of
the points $p_i=(1:1:1:\alpha_i)$ $i=1,2,3,4$ with $\alpha_i$ such
that $4\alpha_i^4+4\lambda\alpha_i+3=0$, of $p_5=(1:\xi_3:\xi_3^2:0)$
and of $p_6=(1:\xi_3^2:\xi_3:0)$ and 4 orbits of points with stabilizer
isomorphic to $\Z/2\Z$, the orbits of the points
$q_1=(\sqrt{\gamma}:\sqrt{\gamma}:1:1)$,
$q_2=(-\sqrt{\gamma}:-\sqrt{\gamma}:1:1)$, $q_3=(\sqrt{\gamma}:
-\sqrt{\gamma}:1:-1)$, $q_4=(-\sqrt{\gamma}:\sqrt{\gamma}:1:-1)$
where $\gamma^2+2\lambda\gamma+1=0$. With respect to the action of
$\mathfrak{S}_4$ the orbit of $p_5$ (resp. $q_3$) coincides with that of $p_6$ (resp. $q_4$).

Let $\pi:X_\lambda^4\ra X_\lambda^4/\mathfrak{S}_4\subset W\mathbb{P}(1,2,3,4)$. Then
$\pi(p_i)=\overline{p}_i$ and $\pi(q_j)=\overline{q}_j$ with
$$\overline{p}_i=(3+\alpha_i:3+3\alpha_i:1+3\alpha_i:\alpha_i), \quad i=1,2,3,4,$$
$$\overline{p}_5=(0:0:1:0),$$
$$\overline{q}_1=(2\sqrt{\gamma}+2:\gamma+4\sqrt{\gamma}+1:2\gamma+2\sqrt{\gamma}:\gamma),$$
$$\overline{q}_2=(-2\sqrt{\gamma}+2:\gamma-4\sqrt{\gamma}+1:2\gamma-2\sqrt{\gamma}:\gamma),$$
$$\overline{q}_3=(0:-\gamma-1:0:\gamma).$$ These are five singular
points ($\overline{p}_i$, $i=1,\ldots, 5$) of type $A_2$  and
three singular points ($\overline{q}_j$, $j=1,2,3$)
of type $A_1$.\\ We observe that the points $\overline{p_5}$ and $\overline{q_3}$ are in the singular
locus of $W\mathbb{P}(1,2,3,4)$, the other singular points are
not. If $\lambda^4=1$, then $X_\lambda^4/\mathfrak{S}_4$ admits
three other singular
points ($(4:6:4:1)$, $(0:-2:0:1)$, $(2i:-2:-2i:1)$), which are the image of the singularities of $X_\lambda^4$.\\
Since the quotient $X_\lambda^4/\mathfrak{A}_4$ is realized as
double cover of $X_\lambda^4/\mathfrak{S}_4$, the singularities on
$X_\lambda^4/\mathfrak{A}_4$ come from the singularities of
$X_\lambda^4/\mathfrak{S}_4$ or from the singularity of the branch
locus $B_\lambda^4$, which is given by
\begin{eqnarray*}
0=det(M_4)&=&
-4e_1^3e_3^3-27e_1^4e_4^2+16e_2^4e_4-4e_2^3e_3^2-27e_3^4+256e_4^3\\ & & +
18e_1^3e_3e_2e_4-80e_1e_3e_2^2e_4-4e_1^2e_2^3e_4 \\ & & +e_1^2e_2^2e_3^2+144e_1^2e_2e_4^2-6e_1^2e_4e_3^2+18e_2e_1e_3^3\\
 & & +144e_2e_4e_3^2-128e_2^2e_4^2-192e_3e_1e_4^2.
\end{eqnarray*}

The singular locus of $B_\lambda^4 \in
W\mathbb{P}(1,2,3,4)$ is given by the point $(0:1:0:0)$ and
$(0:1:0:1/4)$ in $W\mathbb{P}(1,2,3,4)$, which are not points of
$X_\lambda^4/\mathfrak{S}_4$ if $\lambda\neq 1$ (if $\lambda=1$,
$(0:1:0:4)\in X_1^4/\mathfrak{S}_4$ and $(0:1:0:0)\notin
X_1^4/\mathfrak{S}_4$). So if $\lambda\neq 1$ the singularities of
$X_\lambda^4/\mathfrak{A}_4$ are mapped to singularities of
$X_\lambda^4/\mathfrak{S}_4$. The points $\overline{p}_i$,
$\overline{q}_j$ in $X_{\lambda}^4/\mathfrak{S}_4$, $i=1,2,3,4$
and $j=1,2$ are in the branch locus of the double cover
$X_\lambda^4/\mathfrak{A}_4\ra X_\lambda^4/\mathfrak{S}_4$, hence
each of them corresponds to a singular point on
$X_\lambda^4/\mathfrak{A}_4$. The points $\overline{p}_5$,
$\overline{q}_3$ in $X_{\lambda}^4/\mathfrak{S}_4$, are not in the
branch locus of the double cover $X_\lambda^4/\mathfrak{A}_4\ra
X_\lambda^4/\mathfrak{S}_4$, hence each of them corresponds to two
singular points on $X_\lambda^4/\mathfrak{A}_4$ (indeed the two
singular points mapped to $\overline{p}_5$ are the image in
$X_\lambda^4/\mathfrak{A}_4$ of the points in the orbit of $p_5\in
X_\lambda^4$ and $p_6\in X_\lambda^4$ under the action of
$\mathfrak{A}_4$ and analogously the two singular points over
$\overline{q}_3$ are the image in $X_\lambda^4/\mathfrak{A}_4$ of
the points in the orbit of $q_3\in X_\lambda^4$ and $q_4\in
X_\lambda^4$ under the action of $\mathfrak{A}_4$). Hence if
$\lambda^4\neq 1$, the singular points of
$X_\lambda^4/\mathfrak{A}_4$ are 6 points of type $A_2$ and 4
points of type $A_1$, as computed in a general setting by
\cite{Xiao}. If $\lambda^4=1$ we have other four singularities of
type $A_1$, indeed the two singular points of
$X_\lambda^4/\mathfrak{S}_4$, $(4:6:4:1)$ and $(0:-2:0:1)$ are in
the branch locus of the double cover
$X_\lambda^4/\mathfrak{A}_4\ra X_{\lambda}^4/\mathfrak{S}_4$, and
the point $(2i:-2:-2i:2)$ is not.
\endproof

The equations of $T_{\lambda}^4$ are given by
$$
w_4^2-(e_1e_2e_3-e_3^2-e_4e_1^2)w_4-27e_1^4e_4^2 + 4608e_1^3e_2e_3e_4^4 - 4e_1^3e_3^3 - 4e_1^2e_2^3e_4 +
$$
$$
    e_1^2e_2^2e_3^2 + 144e_1^2e_2e_4^2 - 6e_1^2e_3^2e_4 - 80e_1e_2^2e_3e_4 +
    18e_1e_2e_3^3 - 192e_1e_3e_4^2 +
    $$
    $$
    16e_2^4e_4 - 4e_2^3e_3^2 - 128e_2^2e_4^2 +  144e_2e_3^2e_4 - 27e_3^4=0
$$
$$
4(1+\lambda)e_4 - (e_1^4-4e_2e_1^2+4e_3e_1+2e_2^2)=0
$$
in $W{\mathbb P}^4(1,2,3,4,6)$. A desingularization is obtained by blowing up at the singular points: by general facts, we still get a K3 surface.

We end this section with the study of the singularities for $n=4$. The hypersurface ${\mathbb P}^4/{\mathfrak A}_{5}$ in $W{\mathbb P}^5(1,2,3,4,5,10)$ is given by the equation
$$
w_5^2-\theta_5w_5 -\Delta_5=0,
$$
where
$$
\theta_5= e_1e_2e_3e_4-e_1^2e_4^2-e_1e_2^2e_5+2e_1e_4e_5+e_2e_3e_5-e_3^2e_4-e_5^2,
$$
$$
\Delta_5=64e_1^5e_5^3 - 48e_1^4e_2e_4e_5^2 - 32e_1^4e_3^2e_5^2 + 36e_1^4e_3e_4^2e_5 -
    7e_1^4e_4^4 +
    $$
    $$
   + 36e_1^3e_2^2e_3e_5^2 - 2e_1^3e_2^2e_4^2e_5 - 20e_1^3e_2e_3^2e_4e_5 + 5e_1^3e_2e_3e_4^3 - 400e_1^3e_2e_5^3 +
$$
$$
    4e_1^3e_3^4e_5 - e_1^3e_3^3e_4^2 + 40e_1^3e_3e_4e_5^2 - 8e_1^3e_4^3e_5 -
    7e_1^2e_2^4e_5^2 + 5e_1^2e_2^3e_3e_4e_5 +
    $$
    $$
    - e_1^2e_2^3e_4^3 - e_1^2e_2^2e_3^3e_5+ 256e_1^2e_2^2e_4e_5^2 + 140e_1^2e_2e_3^2e_5^2 - 187e_1^2e_2e_3e_4^2e_5 +
    $$
  $$
    36e_1^2e_2e_4^4 + 6e_1^2e_3^3e_4e_5 - 2e_1^2e_3^2e_4^3 + 500e_1^2e_3e_5^3 -
    14e_1^2e_4^2e_5^2 - 157e_1e_2^3e_3e_5^2 +
 $$
 $$
    6e_1e_2^3e_4^2e_5 + 88e_1e_2^2e_3^2e_4e_5 - 20e_1e_2^2e_3e_4^3 + 562e_1e_2^2e_5^3 -
    18e_1e_2e_3^4e_5 + 5e_1e_2e_3^3e_4^2+
    $$
    $$
     - 513e_1e_2e_3e_4e_5^2 + 40e_1e_2e_4^3e_5 - 225e_1e_3^3e_5^2 + 256e_1e_3^2e_4^2e_5 - 48e_1e_3e_4^4 -
  $$
  $$
    624e_1e_4e_5^3 + 27e_2^5e_5^2 - 18e_2^4e_3e_4e_5 + 4e_2^4e_4^3 +
    4e_2^3e_3^3e_5 - e_2^3e_3^2e_4^2 - 225e_2^3e_4e_5^2 +
    $$
    $$
     + 206e_2^2e_3^2e_5^2 +140e_2^2e_3e_4^2e_5 - 32e_2^2e_4^4 - 157e_2e_3^3e_4e_5 + 36e_2e_3^2e_4^3 -
    937e_2e_3e_5^3
    $$
    $$
    + 500e_2e_4^2e_5^2 + 27e_3^5e_5 - 7e_3^4e_4^2 + 562e_3^2e_4e_5^2 - 400e_3e_4^3e_5 + 64e_4^5 + 781e_5^4
$$

By direct inspection, this hypersurface is quasi-smooth and well-formed in weighted projective space. The singularities are given by the intersection of the singularities of $W{\mathbb P}^5(1,2,3,4,5,10)$ and ${\mathbb P}^4/{\mathfrak A}_{5}$. These are given by the curve
$$
w_5^2=4e_4^3(e_2^4-8e_2^2e_4+16e_4^2)
$$
in the weighted projective plane $W{\mathbb P}^2(2,4,10)$ given by $$e_1=e_3=e_5=0,$$
 the point $[0,0,1,0,0,0]$ and two points which satisfy the equation $w_5^2+e^2_5w_5-781e_5^4=0$ in the weighted projective line $W{\mathbb P}^1(5,10)$ defined by the equations $$e_1=e_2=e_3=e_4=0.$$

 These singular loci give rise to points in ${\mathbb P}^4$ which are stabilized by the corresponding isotropy groups. All these points belong to $X^4_{\lambda}$ for generic $\lambda$ except the points with isotropy group the cyclic group of order $5$. Therefore, the singularities of $T_{\lambda}^5$ are given by the point $[0,0,1,0,0,0]$ and the curve
$$
e_1=e_3=e_5=0, \quad  w_5^2=4e_4^3(e_2^4-8e_2^2e_4+16e_4^2).
$$

Since ${\mathfrak A}_5$ acts on $X_{\lambda}^5$ preserving the period, $T_{\lambda}^5$ admits a crepant resolution which is a smooth Calabi-Yau manifold.  The singular point with isotropy group the cyclic group of order three can be resolved via toric crepant resolution. The singular curve is locally ${\mathbb C} \times {\mathbb C}^2 /\left({\Z/2\Z}\right)$.

By standard results on weighted projective space (see, for instance, \cite{Joy}, p. 134), the isotropy group of a point $p$ is cyclic of order $k$, where $k$ is the highest common factor of the weights corresponding to non zero coordinates of $p$. Now, we prove that the quotient $T_{\lambda}^{n+1}$ does not admit a crepant resolution for $n \geq 5$. According to Proposition 6.4.4 in \cite{Joy},  the quotient $T_{\lambda}^{n+1}$  does not have a crepant resolution if it has terminal singularities. We recall that a quotient singularity is terminal if the age is greater than 1. Locally, $T_{\lambda}^{n+1}$ is isomorphic to ${\mathbb C}^{n-1}/G$, where $G$ is a cyclic group of order $k \geq 2$. Let $g$ be a generator of $G$. The age of $g$ is computed as follows.

\begin{lemma}
\label{ageterm}
Let $k \geq 2$. If $n+1=sk+t$ for $0 \leq t \leq k-1$, then the age of $g$ is equal to
$$
\frac{1}{k}\left(\left\lfloor \frac{n+1}{k} \right\rfloor \frac{k(k-1)}{2} +ak\right),
$$
where $a$ is the quotient of the division of $\frac{t(t+1)}{2}$ by $k$.
\end{lemma}

\proof First we recall that every automorphism $\alpha$ of order $r$ of an $m$-dimensional variety linearizes near a component of the fixed locus to a diagonal matrix $diag(\xi_r^{a_1},\xi_r^{a_2},\ldots, \xi_r^{a_m})$, $0\leq a_i< r$ and the age of $\alpha$ is by definition
$$age(\alpha):=\sum_{i=1}^m \frac{a_i}{r}.$$

In our case, by definition, we have $e_k \neq 0$. Consider the subset $U_k$, where $e_k \neq 0$. This is isomorphic to ${\mathbb C}^{n+1}/G$, where $G$ is the cyclic group of order $k$ and the action of $G$ on  ${\mathbb C}^{n+1}$ is given by
\begin{equation}
\label{weights}
(z_{i_1}, \ldots, z_{i_{n+1}}) \rightarrow (\xi_k^{a_{i_1}}z_{i_1}, \ldots, \xi_k^{a_{i_{n+1}}}z_{n+1}),
\end{equation}
where $\xi_k$ is a $k$-th primitive root of unity and $$\{i_1, \ldots, i_{n+1}\}=\left\{1, \ldots, \widehat{k}, \ldots, n,n+1, \frac{n(n+1)}{2}\right\}.$$

Note that $k$ is different from $\frac{n(n+1)}{2}$. If it were equal to it, the point $p=[0,0,\ldots, 0, 1]$ would be singular on $T_{\lambda}^{n+1}$. Clearly, $p$ does not satisfy the equations defining $T_{\lambda}^{n+1}$. Moreover, we can assume that $k$ is different from $n$ and $n+1$. Suppose, for instance, that $k=n$. By definition of $k$, the only non zero coordinates of points with isotropy group ${\mathbb Z}/n{\mathbb Z}$ are $e_n$ and $w_{n+1}$ for odd $n$. In this case, we claim that $w_{n+1}$ is different from zero. If it were, we would have to deal with the point $e_j=0$ for any $1 \leq j \leq n+1$ and $j \neq n$ and $w_{n+1}=0$. It is easy to show that this point does not satisfy the first of the two equations \eqref{eqs}.

Second, we have $a_{i_j} \equiv j$ mod $k$ for $j \neq n+1$, and  $a_{i_{n+1}} \equiv \frac{n(n+1)}{2}$ mod $k$.

As proved in Lemma \ref{Torbifold}, $e_n$ and $e_{n+1}$ can be expressed in terms of the other variables. The projection map $\pi$
$$
W{\mathbb P}^{n+1}\left(1, \ldots, n+1, \frac{n(n+1)}{2}\right) \rightarrow W{\mathbb P}^{n-1}\left(1, \ldots, \widehat{n}, \widehat{n+1},  \frac{n(n+1)}{2}\right),
$$
forgets $e_n$ and $e_{n+1}$.  The equations in \eqref{eqs} have degree $n(n+1)$ and $n+1$, respectively. This implies that  the restriction of $\pi$ to $T_{\lambda}^{n+1}$ is a degree $n+1$ covering. The image of $U_k \cap T_{\lambda}^{n+1}$ under this restriction is isomorphic to ${\mathbb C}^{n-1}/G$.  The weights $b_{r_j}$ of the action of $G$ on ${\mathbb C}^{n-1}$ are given by those of the action \eqref{weights}. More precisely, we have
\begin{equation}
\label{ageI}
\sum_{j=1} b_{r_j}= \sum_{j=1} a_{i_j}-[n+1]_k-[n(n+1)]_k,
\end{equation}
where $[R]_k$ denotes the representative of the class of $R$ modulo $k$ in between $0$ and $k-1$.
The contribution $[n+1]_k$ is due to the removal of the variable $e_{n+1}$, and the contribution $[n(n+1)]_k$ is due to the removal of the variable $e_n$: notice that the first equation in \eqref{eqs} has degree $n+1$ with respect to $e_n$ and the second equation in \eqref{eqs} is linear with respect to the variable $e_{n+1}$. The sum
$$
\sum_j a_{i_j}
$$
can be computed as follows. Write it as
$$
\sum_j^n a_{i_j} +a_{i_{n+1}}= \sum_j a_{i_j} +\left[\frac{n(n+1)}{2}\right]_k,
$$
and note that
$$
\sum_{j=1}^{n}a_{i_j} = \left\lfloor \frac{n+1}{k} \right\rfloor \frac{k(k-1)}{2} +\frac{t(t+1)}{2}.
$$
Therefore, we have
\begin{eqnarray}
\sum_j b_{r_j}  &=&  \left\lfloor \frac{n+1}{k} \right\rfloor \frac{k(k-1)}{2} +\frac{t(t+1)}{2} -\left[\frac{(n+1)(n+2)}{2}\right]_k.
\end{eqnarray}

Since $1+ 2+ \ldots + (n+1) = \frac{(n+1)(n+2)}{2}$, it is an easy exercise to show that
$$
\left[\frac{(n+1)(n+2)}{2}\right]_k=\left[\frac{t(t+1)}{2}\right]_k.
$$

This proves the statement of the lemma. Finally, note that $\sum_j b_{r_j}$ is a multiple of $k$: in fact, we have
\begin{eqnarray}
\sum_j b_{r_j} &\equiv& \sum_j a_{i_j} + \frac{n(n+1)}{2} - (n+1) -n(n+1) \\
&\equiv& \sum _j ^{n+1} j - \frac{(n+1)(n+2)}{2} \\
&\equiv & 0 \, \, \hbox{mod} \, \, k.
\end{eqnarray}

\endproof

\begin{proposition}
\label{noncanonical}
The CY orbifold $T_{\lambda}^{n+1}$ has a crepant resolution if and only if $n=2,3,4$.
\end{proposition}

\proof If $n=2$, then $T_{\lambda}^{n+1}$ is smooth. As proved before, if $n=3,4$ then $T_{\lambda}^{n+1}$ has a crepant resolution.
Assume, now, that $T_{\lambda}^{n+1}$ has a crepant resolution. Thus, it must not have terminal singularities. Suppose, now, that $n \geq 5$. If we show that $T_{\lambda}^{n+1}$ has terminal singularities, it does not have a crepant resolution. For $n \geq 5$ write $n+1=3s+t$, where $0 \leq t <3$. Fix $k=3$. By Lemma \ref{ageterm}, the age of a generator of the cyclic group of order three is given by
$$
age(g)=
\left\{
\begin{array}{c}
s, \quad t=0,1,\\
s+1, \quad t=2.
\end{array}
\right.
$$

An easy calculation shows that the dimension of the fixed locus is $s-2$. Hence, the age of $g^{-1}$ is given by
$$
age(g^{-1})=n-1-dim(Fix(g)) - age(g)=
\left\{
\begin{array}{l}
s +1, \quad t=0,\\
s+2, \quad t=1,2,
\end{array}
\right.
$$
which never equals $1$. In other words, there exists terminal singularities with isotropy group $\Z/3\Z$ for any $n \geq 5$. \endproof

\subsection{The Quotient $Y_{\lambda}^{n+1}$ by the Group $H_n$}\label{sec10}\label{subsec: Y}

By Proposition \ref{prop: Z/(n+1)Z^(n-1) group of autom}, the
group $H_n\simeq (\Z/(n+1)\Z)^{n-1}$ acts on
$X_{\lambda}^{n+1}$ by preserving the period. In this section we take the quotient of
$X_{\lambda}^{n+1}$ by this group. The main point is that the quotient variety is the mirror (as a toric variety) of the family of hypersurfaces of degree $n+1$ in $\mathbb{P}^n$.

For $\lambda =0$ the invariants of the action of $H_n$ are given
by $y_i=x_i^{n+1}$ for $i=1, \ldots, n+1$ and by
$y_0=\prod_{i=1}^{n+1} x_i$. This implies that the quotient of
$X_0^{n+1}$ by $H_n$ is cut out by the following equations in
${\mathbb P}^{n+1}$:
$$
\sum_{i=1}^{n+1} y_i =0, \quad y_0^{n+1} = \prod_{i=1}^{n+1} y_i.
$$

If $\lambda \neq 0$, we have the equations:
$$
\sum_{i=1}^{n+1} y_i =(n+1)\lambda y_0, \quad y_0^{n+1}= \prod_{i=1}^{n+1} y_i.
$$

If we solve for $y_0$ in the first of these equations and substitute in the second one, we get the following equation:
$$
\left(\sum_{i=1}^{n+1} y_i\right)^{n+1}=(n+1)^{n+1}\lambda^{n+1}\prod_{i=1}^{n+1}y_i.
$$

The quotient $Y_{\lambda}^{n+1}:=X_{\lambda}^{n+1}/H_n$ is quite
singular. In fact, it is singular at each point $p$ where the
stabilizer of $p$ in $H_n$ is non-trivial. Notice that a point in
${\mathbb P}^n$ has a non-trivial stabilizer in $H_n$ if at least
two of the coordinates are zero. Let $I$ be a multi-index of
length $k$ where $ 2\leq k \leq n-1$ and denote by $C_{I}$ the
linear subspace $C_I:=\{x_{i_1}=\ldots = x_{i_k}=0\}$. Then the
points of $C_I \cap X_{\lambda}^{n+1}$ have stabilizer isomorphic
to $\left( \Z/ (n+1) \Z\right)^{k-1}$. Moreover the dimension of
$C_I \cap X_{\lambda}^{n+1}$ is equal to $n-k-1$.

\begin{proposition}\label{prop: Y is singular CY}
The quotient $Y_{\lambda}^{n+1}:=X_{\lambda}^{n+1}/H_n$ is a
Calabi-Yau orbifold which admits a projective crepant resolution.
\end{proposition}

\proof Clearly, $Y_{\lambda}^{n+1}$ is an orbifold. Since the elements of $H_n$ preserve the period of $X_{\lambda}^{n+1}$, the singularities of the quotient $Y_{\lambda}^{n+1}$ are Gorenstein, and accordingly, canonical. The canonical sheaf is trivial by the adjunction formula, and the Hodge numbers $h^{0,j}$ are easily seen to be zero - see, for instance, \cite{CK}. As proved in \cite{dais}, Application 5.5, the hypersurface $\Sigma$ given by
$$
y_0^{n+1}= \prod_{i=1}^{n+1} y_i,
$$
is isomorphic to a quotient singularity that admits a projective crepant resolution $\pi: Z_c \rightarrow \Sigma$ in any dimension. The quotient $Y_{\lambda}^{n+1}$ is a hyperplane section of $\Sigma$. Take the divisor $D$ in $Z_c$ given by  $\overline{\pi^{-1}\left(Y_{\lambda}^{n+1}\setminus Sing((Y_{\lambda}^{n+1}\right)}$. Since $\pi$ gives a resolution $\Sigma$, the divisor $D$ yields a resolution of  $Y_{\lambda}^{n+1}$. Moreover, denote by $i$ the inclusion of $D$ in $Z_c$ and $j$ the inclusion of  $Y_{\lambda}^{n+1}$ in $\Sigma$.  By applying twice the adjunction formula, the following holds:
$$
K_D=i^*(K_{Z_c}+D)=i^*\pi^*(K_{\Sigma}+Y_{\lambda}^{n+1})=\pi^*K_{Y_{\lambda}^{n+1}}.
$$

This shows that $Y_{\lambda}^{n+1}$ admits a projective crepant resolution.
\endproof
\begin{remark}\label{rem: Y mirror of X} The variety $Y_\lambda^{n+1}$ is the mirror as a toric variety of $X_\lambda^{n+1}$: see \cite{Bat}. This implies that the Hodge diamond of $Y_\lambda^{n+1}$ is the mirror of that of $X_\lambda^{n+1}$, i.e. for $h^{i,j}(Y_{\lambda}^{n+1})=h^{n-j,i}(X_\lambda^{n+1})$.\\
In case $n=3$, $X_\lambda^{4}$ and $Y_\lambda^{4}$ are K3 surfaces and it is proved  in \cite{dolgachev} that they are mirror also from a lattice theoretic point of view.\\
In case $n=4$, $X_\lambda^{5}$ and $Y_\lambda^{5}$ are Calabi--Yau threefolds and they are mirror in all known definitions of mirror symmetry!
\end{remark}

For $k \geq 3$ and generic $\lambda$, let ${\mathcal I}_{\lambda}^k$ be the variety:
\begin{equation*}
\left\{\left([y_1, \ldots, y_k],[a,b]\right) \in {\mathbb P}^{k-1} \times {\mathbb P}^1 : \right.
\end{equation*}
\begin{equation}
\label{fibration}
\left. \left(y_1+\ldots + y_k\right)^k(b+a)^{k+1} -(k+1)^{k+1}\lambda^{k+1}ab^ky_1 \ldots y_k=0\right\}.
\end{equation}

Denote by $\alpha_k$ the second projection onto ${\mathbb P}^1$.

\begin{theorem}
\label{ikappa}
The variety ${\mathcal I}_{\lambda}^k$ is birational to $Y_{\lambda}^{k+1}$.
\end{theorem}

\proof  Denote by $F^{0}_k$ the fiber over the point $[1,0]$, and let $I_{\lambda}^k$ be the open set obtained from ${\mathcal I}_k$ by removing $F^{0}_k$. Consider the rational map
$$
\beta_{k+1}: Y_{\lambda}^{k+1} \dashrightarrow {\mathbb P}^1,
$$
which maps $(y_1, \ldots, y_{k+1})$ to $(y_{k+1}, y_1 + \ldots + y_k)$. Let $V_{k+1}$ the union of $\beta_{k+1}^{-1}([1,0])$ and $\Pi_{k+1}$, where $\Pi_{k+1}$ is the codimension two subspace given by
$$
y_1+\ldots+y_k=0, \quad y_{k+1}=0.
$$
Set ${\mathcal Y}_{\lambda}^{k+1} = Y_{\lambda}^{k+1} \setminus V_{k+1}$. The restriction of $\beta_{k+1}$ to ${\mathcal Y}_{\lambda}^{k+1}$ is well-defined because we removed the indeterminacy locus. The morphism between ${\mathcal Y}_{\lambda}^{k+1}$ and  $I_{\lambda}^k$  is given by
$$
\left[y_1, \ldots, y_k, y_{k+1}\right] \rightarrow \left(\left[y_1, \ldots, y_k],[y_{k+1},y_1+\ldots + y_k\right]\right).
$$
Its inverse is given by
$$
\left(\left[y_1, \ldots, y_k\right],\left[a,b\right]\right) \rightarrow \left[y_1, \ldots, y_k, \frac{a}{b}(y_1+\ldots + y_k)\right].
$$
\endproof

\begin{remark} We recall that $Y_{\lambda}^3$ is a quotient of $X_{\lambda}^3$. Theorem \ref{ikappa} shows that we can reconstruct $Y_{\lambda}^{k+1}$ as a fibration over ${\mathbb P}^1$
with general fiber $Y_{\mu}^k$ for suitable $\mu$, which can be computed explicitly. Thus, all $Y_{\lambda}^{k+1}$ originate from $X_{\lambda}^3$.
\end{remark}

In order to get more information on the geometry of $Y_{\lambda}^{n+1}$, we try to understand whether they may be viewed as coverings of projective space.

Let us consider the Calabi-Yau $Y_{\lambda}^{n+1}$ of dimension $n-1$
with equation:
$$\left(\sum_{i=1}^{n+1}y_i\right)^{n+1}-\gamma_n\prod_{i=1}^{n+1}y_i=0,$$
where $\gamma_n= ((n+1)\lambda)^{n+1}$.

Let us consider the birational map
$\pi:\mathbb{P}^{{n}}\ra\mathbb{P}^{{n-1}}$ given by
$$\pi: (y_1:
\ldots: y_{n+1}) \mapsto (y_1+y_2:y_3:y_4:\ldots:y_{n+1}).$$ Set
$$
y_i=\mu_{i-2}(y_1+y_2),\ i=3,\ldots {n+1}
$$
where $(1:\mu_1:\mu_2:\ldots:\mu_{{n-1}})$ are
coordinates on ${\mathbb P}^{{n-1}}$.

If we substitute these equations in the pencil and put $y_2=1$, we
get
$$
Ay_1^2+y_1(2A-\gamma_n B)+A=0,
$$
where $A=(1+\sum_{i=1}^{{n-1}}\mu_i)^{n+1}$ and
$B=\prod_{i=1}^{{n-1}}\mu_i$. This shows $Y_{\lambda}^{n+1}$
as a double cover of $\mathbb{P}^{{n-1}}$ branched along the
discriminant of the equation $Ax_1^2+x_1(2A-\gamma_n B)+A=0$ with variable
$x_1$. The discriminant (in the homogeneous coordinates
$(\mu_0:\mu_1:\ldots :\mu_{{n-1}})$) is given by
$$
\gamma_n \prod_{i=1}^{{n-1}}\mu_i\left(\gamma_n \mu_0^2\prod_{i=1}^{{n-1}}\mu_i-4\left(\sum_{i=0}^{{n-1}}\mu_i\right)^{n+1}\right).$$
Hence we obtain a model of $Y_{\lambda}^{n+1}$ as a double cover of
$\mathbb{P}^{{n-1}}$ with branch locus of degree $2n$: \begin{equation}\label{equation: double cover}z^2=\gamma_n \prod_{i=1}^{{n-1}}\mu_i\left(\gamma_n\mu_0^2\prod_{i=1}^{{n-1}}\mu_i-4\left(\sum_{i=0}^{{n-1}}\mu_i\right)^{n+1}\right).\end{equation} The branch
locus is the union of ${n-1}$ hyperplanes $\mu_i=0$ for
$i=1,\ldots {n-1}$ and an irreducible singular hypersurface of
degree ${n+1}$, singular along the $n-3$ dimensional linear space
$$L_{n-3}:=\left\{\begin{array}{l}\mu_0=0\\\sum_{i=1}^{{n-1}}\mu_i=0\end{array}\right.,$$
and along the $\left(\begin{array}{c} {n-1}\\2\end{array}\right)$
linear spaces of dimension $n-4$, namely:
$$R_{h,k}:= \left\{\begin{array}{l}\mu_h=0,\\\mu_k=0,\\ \sum_{i=0}^{{n-1}}\mu_i=0,\end{array}\right.$$
where $h,k$ are chosen so that $h< k$,
$h,k\in\{1,\ldots,{n-1}\}$.\\

Let us consider the pencil of hyperplanes containing
$L_{n-3}$: $\tau \sum_{i=1}^{{n-1}} \mu_i=\mu_0$. Each
hyperplane $H$ of this pencil intersects the branch locus in
$L_{{n-3}}$ and in a hypersurface in $H$ of degree $2n-2$. If we
take into account the pencil of hyperplanes through $L_{n-3}$, we
obtain a fibration over $\mathbb{P}^1_{\tau}$ with fibers the varieties
$V_{n-2}$ of dimension $n-2$. These varieties are double cover of
$H(=\mathbb{P}^{n-2})$ branched along a hypersurface of degree
$2n-2$. The equation of such a fibration is
$$F_1:=\left\{z^2=\prod_{i=1}^{{n-1}}\mu_i\left[\gamma_n \tau^2\prod_{i=1}^{{n-1}}\mu_i-4(\tau+1)^{n+1}(\sum_{i=1}^{{n-1}}\mu_i)^{{n-1}}\right]\right\}.$$
The generic fiber is a variety which is a double cover of $\mathbb{P}^{n-2}$ branched along the union of $n-1$ hyperplanes and of the hypersurface $Y^{n-1}_{\nu(\tau, \lambda)}$ where $$\nu(\tau, \lambda)=\frac{\tau^2\gamma_{n}}{4(\tau+1)^{n+1}}=\frac{\tau^2}{4}\left(\frac{(n+1)\lambda}{(\tau+1)}\right)^{n+1}.$$

Analogously, one can construct the fibration $F_2$ associated to
the pencil of hyperplanes ($\tau \mu_i=\mu_j$, $i\neq j$,
$i,j=1,\ldots {n-1}$) through the ${n-3}$ linear subspaces of $M_{h,k}$, which are
the intersection of the hyperplanes $\mu_h=0$ and $\mu_k=0$. For $i=2$ and
$j=1$, we obtain
$$F_2:=\left\{z^2=\tau\prod_{i=3}^{{n-1}}\mu_i\left[\gamma_n \tau \mu_0^2\mu_2^2\prod_{i=3}^{{n-1}}\mu_i-4\left(\mu_0+\mu_2(1+\tau)+\sum_{i=3}^{{n-1}}\mu_i\right)^{n+1}\right]\right\}.
$$

\subsubsection{The Curve $Y_{\lambda}^3$}
The group $H_2$ acts on the elliptic curve $X_\lambda^3$ fixing the period and thus is generated by a translation by a point of order 3 on $X_\lambda^3$. The action of $H_2$ is fixed point free and the quotient curve $Y_\lambda^3$ is the elliptic curve given by the equation $(y_1+y_2+y_3)^3-3\lambda y_1y_2y_3=0$. We will analyze better the properties of $Y_\lambda^3$ and its relation with $X_\lambda^3$ in the Section \ref{sec: curve quotient by the semidirect product}.

\subsubsection{The Surface $Y_{\lambda}^4$}\label{subsub K3 Y}

The surface $Y_{\lambda}^4$ - denoted $Y_\lambda$ - is singular. By definition, it is obtained as a  quotient of a K3 surface ($X_\lambda^4$) by a group of automorphisms ($H_3$) acting trivially on the period. The automorphisms of K3 surfaces, which are trivial on the period, are called {\it symplectic automorphisms} - indeed they preserve the symplectic structure of the surface - and are well-known and studied in the literature. In particular, the finite groups acting symplectically on a K3 surface are classified (cf. \cite{Nikulin symplectic}, \cite{Mu}, \cite{Xiao}). Let $G$ be a finite group of automorphism on a K3 surface $S$. The desingularization $\widetilde{S/G}$ of $S/G$ is a K3 surface if and only if $G$ is a group of symplectic automorphisms.  If $G$ is symplectic, let us consider the minimal primitive sublattice of $NS(\widetilde{S/G})$ which contains the curves arising by the desingularization of $S/G$. It depends on $G$, but not on $S$ for almost all the groups acting symplectically on a K3 surface (the unique exceptions are $Q_8$ and $T_{24}$), and in this case we will denote it by $M_G$. In \cite{Xiao}, the lattices $M_G$ are computed for all admissible groups $G$.\\

Since the automorphisms of $H_3$ on $X^4_{\lambda}$ are symplectic, the desingularization of
$Y_{\lambda}$ is a K3 surface.\\
The six points of $Sing(Y_{\lambda})$ are
$$\begin{array}{ccc}p_1=(0:0:-1:1),& p_2=(0:-1:0:1),&
p_3=(-1:0:0:1),\\
p_4=(0:-1:1:0),& p_5=(-1:0:1:0), & p_6=(-1:1:0:0).\end{array}$$
These singularities are of type $A_3$ (because they are the image
of the points $q_i$ on $X^4_{\lambda}$ which are stabilized by a
subgroup isomorphic to $\Z/4\Z$ in $H_3$).\\
Since $Y_\lambda$ is the desingularization of the quotient of a K3
surface with the group $(\Z/4\Z)^2$ of symplectic automorphisms, we have
$M_{(\Z/4\Z)^2}\hookrightarrow NS(Y_\lambda)$.
Moreover, $Y_\lambda$ has a singular model as a quartic in
$\mathbb{P}^3$ and the singularities are obtained via
contraction of the curves in $M_{(\Z/4\Z)^2}$; hence $\Z L\oplus
M_{(\Z/4\Z)^2}\hookrightarrow NS(Y_{\lambda})$
with finite index, where $L$ is a pseudo-ample divisor with $L^2=4$.\\
This index can be $1,2,4$. In fact, the lattice $\Z L\oplus
M_{(\Z/4\Z)^2}$ is the N\'eron--Severi group of a certain K3
surface. As proved in \cite{elliptic io}, the
lattice $U\oplus M_{(\Z/4\Z)^2}$ is the N\'eron--Severi group of
an elliptic K3 surface with $6I_4$ as reducible fibers . Let us
consider the sublattice $\langle 4\rangle\oplus
M_{(\Z/4\Z)^2}\simeq \left(\begin{array}{c}1\\2\end{array}\right) \oplus
M_{(\Z/4\Z)^2}$ of $U\oplus M_{(\Z/4\Z)}^2$. It is a primitive
sublattice of $U\oplus M_{(\Z/4\Z)^2}$ and hence it is a primitive
sublattice of $\Lambda_{K3}$. Thus, it appears as the N\'eron--Severi
group of a K3 surface.\\
One can construct lattices $R$ such that $\Z L\oplus
M_{(\Z/4\Z)^2}$ has index 2 or 4 in $R$. These lattices are
primitively embedded in $\Lambda_{K3}$ because their length
($l(R)$) is 1 ($=l(M_{(\Z/4\Z)^2}+1-2$) and their rank is rank$R$=19. Hence
they satisfy the inequality $l(R)\leq
rank(\Lambda_{K3})-rank(R)-2$ which guarantees that $R$ is
primitively embedded in $\Lambda_{K3}$ (cf. \cite{Nikulin bilinear}, Theorem
1.14.4).\\

We recall that $Y_{\lambda}$ is the mirror (in the sense of the $L$-polarized K3 surfaces) of the quartic hypersurface in $\mathbb{P}^3$ (cf. \cite{dolgachev} and Remark \ref{rem: Y mirror of X}). This implies that $T_{Y_\lambda}\simeq \langle 4\rangle\oplus U$. Hence $NS(Y_{\lambda})^{\vee}/NS(Y_{\lambda})=\Z/4\Z$ and
$NS(Y_\lambda)\hookleftarrow \Z L\oplus M_{(\Z/4\Z)^2}$ with index 4.\\
We recall that $M_{(\Z/4\Z)^2}$ is generated by the classes
$M_{j,i}$, the standard basis of the $j$-th copy of the lattice $A_3$ and by
the classes
$$d_1=f_1+f_2+f_3+f_4+2f_5,$$
$$d_2=f_1+2f_2+3f_3+f_5+f_6,$$
where $f_h=(M_{h,1}+2M_{h,2}+3M_{h,3})/4.$\\
Let us consider the class
$$
v=f_1-f_2+f_3=
$$
$$
=\frac{1}{4} \left(M_{1,1}+2M_{1,2}+3M_{1,3}-M_{2,1}-2M_{2,2}+\right.
$$
$$
\left. - 3M_{2,3}+M_{3,1}+2M_{3,2}+3M_{3,3}\right).
$$
Let us denote by $T$ the lattice spanned by the generators of ${\mathbb Z}L \oplus M_{(\Z/4\Z)^2}$ and the class $L/4+v$.

\begin{theorem} The overlattice $T$ is unique up to isometry. In particular, the generators of ${\mathbb Z}L \oplus M_{(\Z/4\Z)^2}$ and $L/4+v$ form a ${\mathbb Z}$-basis of $NS(Y_{\lambda})$ for generic $\lambda$.
\end{theorem}

\proof The lattice $T$ is an even overlattice of ${\mathbb Z}L \oplus M_{(\Z/4\Z)^2}$ of index four. Its signature is $(1,18)$. The uniqueness follows from \cite{Nikulin bilinear}, Corollary 1.13.3.
\endproof

\begin{remark} The Picard number of $X_\lambda^4$ and $Y_\lambda^4$ is the same, since there exists a rational dominant map from $X^4_\lambda$ to $Y_\lambda^4$. For certain specific values of $\lambda$, the Picard number of $Y_\lambda$ and $X^4_{\lambda}$ is 20. In particular, if $\lambda^4=1$, $X^4_\lambda$ is a singular quartic and the Picard group of its desingularization contains an extra class (the exceptional curve of the blow-up of the singular points). Thus, if $\lambda^4=1$, the Picard number of $Y_\lambda^4$ jumps to 20. Similarly, if $\lambda=0$ , the variety $X_0^4$ has Picard number 20 (as proved before) and thus the Picard number of $Y_0^4$ jumps to 20.\end{remark}

\begin{remark} The surfaces $Y_\lambda$ for a generic $\lambda$ are not Kummer surfaces; indeed there exists no even lattice $L$ such that $T_{Y_\lambda}\simeq L(2)$.\end{remark}

The equation \eqref{equation: double cover} exhibits $Y_{\lambda}$ as a double cover of $\mathbb{P}^2$ and the branch locus is the union of two lines and a quartic curve, namely:
$$
l:=\mu_1=0, \,\,\, m:=\mu_2=0 \,\,\, Q:=[2(\mu_0 +\mu_1
+\mu_2)^4- \gamma_3 \mu_1 \mu_2\mu_0^2].$$

The branch locus
has four singular points: $$p_1:=(0:0:1),\ \ p_2=(0:1:-1),\ \
p_3=(1:0:-1),\ \ p_4=(1:-1:0)
$$ where $p_1= l\cap m$, $p_2=Q\cap l$, $p_3=Q\cap m$, and $p_4$
is a singular point on the quartic $Q$. \\
To resolve the singularity $p_1$, it suffices to blow up once the
point $p_1$, thus introducing an exceptional divisor. To resolve the
singularity of the quartic, it suffices to perform a blow-up. The
exceptional divisor of this blow-up is tangent to the strict
transform of the quartic. It has multiplicity 2 and hence it is
not in the branch locus. To separate the strict transform of the
quartic from that of a line ($l$ or $m$), it is necessary to blow-
up 4 times, introducing four exceptional divisors for each line.
Hence, to resolve the singularities of the sextic (and so of its
double cover), we introduce ten exceptional divisors. In particular,
the double cover of the plane branched on the sextic $l\cup m\cup
Q$ has two singularities of type $A_1$ and two singularities of
type $A_4$.\\

Now, we consider elliptic fibrations on $Y_{\lambda}$ which specialize to $F_1$ and $F_2$ in case
$n=3$.\\
The fibration $F_1$ specializes to the
elliptic fibration
$$z^2=\mu_1\left[\gamma_3 \tau^2\mu_1-4(\tau+1)^4(\mu_1+1)^{2}\right],$$
where we put $\mu_2=1$ and $\tau$ is a parameter on
$\mathbb{P}^1$. After elementary transformations, an equation of such an elliptic fibration is given by
$$z^2=\mu_1\left[\mu_1^2-\mu_1\frac{1}{4}(\gamma_n \tau^2(1-\tau)^2-8)+1\right],$$
which exhibits this elliptic fibration as a specialization of the
elliptic fibration analyzed in \cite[Section 4.6]{bert Nikulin involutions}.  In
particular, this fibration has one fiber of type $I_{16}$ on the
point at infinity, two fibers of type $I_2$ (on $\tau=0$ and
$\tau=1$) and four fibers of type $I_1$. It admits a 4-torsion
section (it is of type given in
\cite{elliptic io}, Proposition 2.2, putting $f=1$,
$e=\frac{1}{2}\sqrt{\gamma_3}\tau(1-\tau)$). Hence the K3 surfaces
$Y_{\lambda}$ admits a symplectic automorphism of order 4. The square of
this automorphism is a Morrison-Nikulin involution which
specializes to the one studied in \cite{bert Nikulin involutions}.

The fibration $F_2$ specializes to the genus
one fibration
$$z^2=\tau\left[\gamma_3 \tau \mu_0^2\mu_2^2-4\left(\mu_0+\mu_2(1+\tau)\right)^4\right].$$

\subsubsection{The Three-fold $Y^5_{\lambda}$}\label{sec: Y5}
For a generic $\lambda$ the 3-fold $Y^5_\lambda$ is a 2:1 cover of $\mathbb{P}^3$ branched along the union of three planes and a quintic surface.\\
The fibration $F_1$ on $Y^5_\lambda$ is a fibration with K3 surfaces as generic fiber. Indeed the fiber of the fibration $F_1$ over the point $\overline{\tau}$ is the double covers of $\mathbb{P}^2$ with coordinates $(\mu_1:\mu_2:\mu_3)$ $$z^2=\mu_1\mu_2\mu_3(\gamma_4\overline{\tau}\mu_1\mu_2\mu_3-4(\overline{\tau}+1)^5(\mu_1+\mu_2+\mu_3)^3).$$
If $\overline{\tau}\neq 0,-1$ this is the double cover of $\mathbb{P}^2$ branched along three lines ($m_1$, $m_2$, $m_3$) and a smooth cubic $C_3$. The lines $m_i$ are inflectional tangents to $C_3$.  The double cover of this singular sextic is of course singular, but it suffices to blow up three times the 3 intersection points between the lines and the cubic and once the 3 intersection points of two lines to obtain a curve such that the double cover branched along this curve is smooth. Hence the singularities of the double cover of the singular sextic are ADE. This implies that the generic fiber is a K3 surface. The fiber over $\overline{\tau}=-1$ is made up of two projective planes and it is not connected, the ones over $\overline{\tau}=0$ and $\overline{\tau}=\infty$ are a double cover of $\mathbb{P}^2$ branched along four lines.

The fibration $F_2$ on $Y^5_\lambda$ is a fibration with K3 surfaces as generic fiber. Indeed the fiber of the fibration $F_2$ over the point $\overline{\tau}$ is the double cover of $\mathbb{P}^2$ with coordinates $(\mu_0:\mu_2:\mu_3)$
$$z^2=\tau\mu_3\left[\gamma_3 \tau \mu_0^2\mu_2^2\mu_3-4\left(\mu_0+\mu_2(1+\tau)+\mu_3\right)^{5}\right].$$
The generic fiber is a double cover of $\mathbb{P}^2$ branched along the union of a singular quintic and a line, which meets the quintic in one point with multiplicity 5. To obtain a smooth model of this double cover of the complex projective plane, we consider the double cover branched along a blow-up of the branch locus. In particular, it suffices to blow up the singular point of the quintic  two times and the intersection points between the line and the quintic five times. This double cover is smooth and hence the generic fiber of this fibration is a K3 surface.

\subsection{The Quotient by ${\mathfrak S}_{n+1} \rtimes H_n$}\label{sec13}\label{subsec: quotients by G}
The full automorphism group of $X_\lambda^{n+1}$ is ${\mathfrak S}_{n+1} \rtimes H_n$ and ${\mathfrak A}_{n+1} \rtimes H_n$ is the subgroup which preserves the period. Here we describe the quotient of $X_\lambda^{n+1}$ by these two groups. In particular we prove:
\begin{proposition} The quotient $X_\lambda^{n+1}/{\mathfrak S}_{n+1} \rtimes H_n$ is birational to the hypersurface $e_1^{n+1}-(n+1)^{n+1}\lambda^{n+1}e_{n+1}=0$ in $W\mathbb{P}(1,2,\ldots, n+1)$ and the quotient $X_\lambda^{n+1}/{\mathfrak A}_{n+1} \rtimes H_n$ is birational to the complete intersection of the hypersurfaces $w_{n+1}^2-w_{n+1}\theta_{n+1}-\Delta_{n+1}=0$ and $e_1^{n+1}-(n+1)^{n+1}\lambda^{n+1}e_{n+1}=0$ in $W{\mathbb P}^n\left(1, 2, 3, \ldots, n+1, \frac{n(n+1)}{2}\right)$. If $n=2,3,4$ there exists a crepant resolution of $X_\lambda^{n+1}/{\mathfrak A}_{n+1} \rtimes H_n$ which is a Calabi--Yau variety.\end{proposition}

Let $q_n: X_{\lambda}^{n+1} \rightarrow Y_{\lambda}^{n+1}$ be the quotient map. Denote by $\pi_n: {\mathbb P}^n \rightarrow {\mathbb P}^n/{\mathfrak S}_{n+1} \cong W{\mathbb P}^n(1, 2, 3, \ldots, n+1)$ the quotient by the symmetric group. Let $D_{\lambda}^{n+1}$ be the image of $X_{\lambda}^{n+1}$ under the composition $\Phi_n=\pi_n \circ q_n$. Clearly, the equation of $D_{\lambda}^{n+1}$ is given by
\begin{equation}
e_1^{n+1}-(n+1)^{n+1}\lambda^{n+1}e_{n+1}=0.
\label{symmtotquot}
\end{equation}
Consider, now, the quotient map $Q_n: X_{\lambda}^{n+1} \rightarrow Q_{\lambda}^{n+1}$ onto the quotient of $X_{\lambda}^{n+1}$ by the semidirect product of ${\mathfrak S}_{n+1}$ and $H_n$.

\begin{lemma}
\label{birequiv}
The maps $\Phi_n$ and $Q_n$ are birationally equivalent.
\end{lemma}

\proof Let $x_1$ and $x_2$ be two points of $X_{\lambda}^{n+1}$ such that $x_1=n_1n_2\cdot x_2$ for $n_1 \in  H_n$ and $n_2 \in {\mathfrak S}_{n+1}$. Then, we have $q_n(x_1)=q_n(n_2\cdot x_2)$. Since $y_i=x_i^{n+1}$, we have $q_n(n_2\cdot x_2)=n_2\cdot q_n(x_2)$. Thus, we obtain
$$
\Phi_n(x_1)=\pi_n \circ q_n (x_1)= \pi_n (n_2 \cdot q_n (x_2)) =   \pi_n \circ  q_n (x_2) = \Phi_n(x_2).
$$

Therefore, there exists a map from $Q_{\lambda}^{n+1}$ to $D_{\lambda}^{n+1}$. Since $\Phi_n$ and $Q_n$ have the same degree, this map must be one-to-one.
\endproof

By Lemma \ref{birequiv}, we can study the quotient by the semidirect product of ${\mathfrak S}_{n+1}$ and $H_n$ up to birationality. This is given by Equation \eqref{symmtotquot}. It is easy to check that $D_{\lambda}^{n+1}$ is quasi-smooth and well-formed. This means that the singularities of $D_{\lambda}^{n+1}$ are those of weighted projective space. Since it has degree $n+1$ in $W{\mathbb P}^n(1, 2, 3, \ldots, n+1)$, it is not Calabi-Yau, i.e., the canonical sheaf is not trivial. As in Lemma \ref{birequiv}, the quotient map of $X_{\lambda}^{n+1}$ onto the semidirect product of ${\mathfrak A}_{n+1}$ and $H_n$ is birational to the map $X_{\lambda}^{n+1} \rightarrow {\mathcal A}_{n+1}$ onto the Calabi-Yau orbifold defined by the equations
$$
w_{n+1}^2-w_{n+1}\theta_{n+1}-\Delta_{n+1}=0, \quad e_1^{n+1}-(n+1)^{n+1}\lambda^{n+1}e_{n+1}=0
$$
in $W{\mathbb P}^n\left(1, 2, 3, \ldots, n+1, \frac{n(n+1)}{2}\right)$. Analogously to the quotient by the alternating group ${\mathfrak A}_{n+1}$, this orbifold admits a crepant resolution for $n=2,3,4$.

\subsubsection{The Case $n=2$}\label{sec: curve quotient by the semidirect product}
By Proposition \ref{prop: An and Sn groups of autom}, the groups $H_2\simeq
\Z/3\Z$ and ${\mathfrak A}_3\simeq \Z/3\Z$ fix the period of the curve $X_{\lambda}^3$. Thus $\mathfrak{A}_3\rtimes H_2\simeq (\Z/3\Z)^2$ is a group of translations on the elliptic curve, and in particular the group of the translation by the 3-torsion of the elliptic curve. Here we consider the quotient $X^3_{\lambda}/({\mathfrak A}_3\rtimes H_2)$, which is isomorphic to the quotient
$Y_{\lambda}^3/{\mathfrak A}_3$.

\begin{proposition}\label{prop: quotient n=2} Let $G\simeq (\Z/3\Z)^2$ be the group of automorphisms of
$\mathbb{P}^2$ generated by $\alpha:(x:y:z)\mapsto (x:\xi_3
y:\xi_3^2z)$ and $\beta:(x:y:z)\mapsto (y:z:x)$. The group $G$ is a group of
automorphisms of the pencil $X^3_{\lambda}$ and the quotient pencil
$X^3_{\lambda}/G$ is isomorphic to $X_{\lambda}^3$.\\
In particular, the quotient of the pencil
$Y_\lambda^3$ by $\beta$ is isomorphic to
the pencil $X^3_{\lambda}$.
\end{proposition}

\proof If $\lambda^3\neq 1$, the curves in the pencil
$X_{\lambda}$ are smooth. They are elliptic curves and the
automorphisms $\alpha$ and $\beta$ have no fixed points. Hence
these automorphisms correspond to a translation by points of order
three (choosing $(1:\xi_6:0)$ as origin of the elliptic curve the automorphism  $\alpha$ is the translation by the point $(1:-1:0)$ and $\beta$ is a translation by the point
$(0:1:\xi_6)$). Hence, for each $\lambda^3\neq 1$ the quotient
$X^3_{\lambda}/G$ is the quotient of an elliptic curve by its
torsion group of order 3. This quotient corresponds to the
multiplication by 3 on the elliptic curve, and hence we obtain an
isomorphic curve.\\
Since $\alpha$ and $\beta$ commute, the quotient
$X_{\lambda}^3\simeq X_{\lambda}^3/G$ is the quotient
$(X_{\lambda}^3/\alpha)/\beta\simeq Y_{\lambda}^3/\beta$.\endproof

We can explicitly describe the isomorphism between the curve $X_{\lambda}^3$ and the curve $X_{\lambda}^3/G$. If we set $x_3=1$ in the equation of $X_{\lambda}^3$, we obtain a pencil of elliptic curves which has the following Weierstrass form:
$$
u^3+\left( - \frac{27}{16}\lambda^4+\frac{27}{2}\lambda\right)u-\frac{27}{32}\lambda^6+\frac{27}{4}-\frac{135}{8}\lambda^3+v^2=0.
$$

The quotient of $X_{\lambda}^3$ by $\alpha$ is given by $(y_1+y_2+y_3)^3-27\lambda^3y_1y_2y_3=0$. The quotient by $\beta$ is given by $$V^2=\frac{4e_1^6}{27\lambda^3}-\frac{2e_1^4e_2}{3l^3}+e_1^2e_2^2-4e_2^3-\frac{e_1^6}{27l^6},$$
where $e_i$ is the $i$-th elementary symmetric functions.  If we set $z=\Delta$ and $x=e_2$, $y=e_1^2$ and, after that, we pass to affine coordinates $s,t$,  we get a pencil of elliptic curves with Weierstrass form:
$$
s^3+(-177147\lambda^{24}+1417176\lambda^{21})s-28697814\lambda^{36}-573956280\lambda^{33}+229582512\lambda^{30}+t^2=0.
$$
The automorphism is thus given by
$$
(x,y) \rightarrow \left( \frac{-\alpha u - 243 \lambda^{12}}{2916\lambda^{12}}, \frac{\alpha^{3/2}v^3}{78732\lambda^{18}}\right),
$$
where
$$
\alpha^2=\frac{-27\lambda^4+216\lambda}{16(-177147\lambda^{24}+1417176\lambda^{21})},
$$
and
$$u=-\frac{3(4y^2+4\lambda xy-4y+4+4\lambda x+\lambda^2x^2)}{4x^2},
 $$
 $$
 v=\frac{9(2y^2+2-y^2x\lambda-2y+2xy\lambda-y\lambda^2x^2+x^3+\lambda x+\lambda^2x^2)}{2x^3}.
$$

\subsubsection{The Case $n=3$}\label{subsub K3 quotient A4} The quotient $X_\lambda^4/\mathfrak{A}_4\rtimes H_3$ is the quotient of a K3 surface by a finite group of symplectic automorphism. In particular, the singularities of this quotient are one singularity of type $D_4$, 6 of type $A_2$ and 2 of type $A_1$ as computed by \cite{Xiao}. The desingularization of $X_\lambda^4/\mathfrak{A}_4\rtimes H_3$ is a K3 surface and its Picard number is 19 if $\lambda$ is generic (the Picard number is the same as that of $X_\lambda^4$). The Picard group is an overlattice of finite index of $\Z L\oplus M_{\mathfrak{A}_4\rtimes H_3}$ for an ample divisor $L$, where $M_{\mathfrak{A}_4\rtimes H_3}$ has discriminant group $(\Z/2\Z)^4\times (\Z/3\Z)^4$ and is an overlattice of index 3 of $D_4\oplus 6A_2\oplus 2A_1$.

\subsubsection{The case $n=4$} The quotient $X_\lambda^5/\mathfrak{A}_5\rtimes H_4$ has a desingularization which is a Calabi-Yau threefold. We will show in Section \ref{quozientone} that the Hodge numbers of any smooth resolution are $h^{1,1}=15$ and $h^{2,1}=5$. This coincides with previous results in \cite{St}, where the Hodge numbers are obtained via representation theory of ${\mathfrak A}_5$. Our method is purely geometric.

\subsection{Automorphisms Not Preserving the Period of $X^{n+1}_{\lambda}$}\label{sec16}\label{subsec: quotient by tau}

Let $\tau$ be the automorphism of $X_{\lambda}^{n+1}$ which does not preserve the period. It is an order two element which generates $\mathfrak{S}_{n+1} \rtimes H_n$ together with the group $\mathfrak{A}_{n+1} \rtimes H_n$, which preserves the period of $X_{\lambda}^{n+1}$. We may assume that $\tau$ is the transposition $(12)$.

If $n=2$, then $X_\lambda^{3}/<\tau>$ is the quotient of an elliptic curve by the hyperelliptic involution, so it is a rational curve.

If $n+1$ is even, the fixed locus of $\tau$ is the divisor $X_\lambda^{n+1}\cup\{x_1=x_2\}$, which is a hypersurface of degree $n+1$ in $\mathbb{P}^{n-1}$; the quotient $X_{\lambda}^{n+1}/<\tau>$ is a smooth variety.

If $n+1$ is odd, the fixed locus of $\tau$ consists of the point $(1:-1:0:\ldots:0)$ and the divisor $X_\lambda^{n+1}\cup\{x_1=x_2\}$, which is a hypersurface of degree $n+1$ in $\mathbb{P}^{n-1}$; the quotient $X_{\lambda}^{n+1}/<\tau>$ is singular. Let $\pi_{(1:-1:0:\ldots:0)}:\mathbb{P}^n\ra \mathbb{P}^{n-1}$ be the projection of $\mathbb{P}^n$ from the point $(1:-1:0:\ldots:0)$ to the $n-1$ dimensional projective space $x_1=x_2$. The restriction of $\pi_{(1:-1:0:\ldots:0)}$ to $X_{\lambda}^{n+1}$ is invariant with respect to the automorphism $\tau$ and exhibits the quotient $X_{\lambda}^{n+1}/<\tau>$ as an $(n/2):1$ covering of $\mathbb{P}^{n-1}$. Indeed, the generic line through $p=(1:-1:0:\ldots:0)$ intersects $X^{n+1}_{\lambda}$ in $n=n+1-1$ points away from $(1:-1:0:\ldots:0)$ and thus $(\pi_{(1:-1:0:\ldots:0)})_{|X_\lambda^{n+1}}$ is generically $n:1$. Since it is invariant under $\tau$, it factorizes through the quotient by $\tau$. Thus $X_\lambda^{n+1}/<\tau>$ is a covering of $\mathbb{P}^{n-1}$. The degree of such a covering depends on the degree of the map $(\pi_{(1:-1:0:\ldots:0)})_{|X_\lambda^{n+1}}$ and of the quotient by $\tau$.

\begin{proposition}
If $n+1$ is even, the quotient $Z:=X_\lambda^{n+1}/<\tau>$ is a smooth Fano variety. Moreover, it is a degree $\frac{n+1}{2}$ covering of $\mathbb{P}^{n-1}$, where the ramification divisor $R$ is linearly equivalent to $(1-n)K_Z$.
\end{proposition}

\proof Let $p: X_\lambda^{n+1} \rightarrow Z$ be the degree two map, which is ramified over the degree $n+1$ hypersurface given by $X_\lambda^{n+1} \cap \{x_1=x_2\}$. By Riemann-Hurwitz's formula, we have $p^*(-K_Z)=h$, where $h$ is the restriction of the hyperplane class of ${\mathbb P}^n$ to $X_\lambda^{n+1}$. Thus, $-K_Z$ is ample, and $Z$ is a smooth Fano variety.

Let $D$ be the divisor on $Z$ associated with the map $u:Z \rightarrow {\mathbb P}^{n-1}$. This means $u^*(H_{n-1})=D$, where $H_{n-1}$ is the hyperplane class on ${\mathbb P}^{n-1}$. If we show that $p^*(D)=h$, then $D$ is linearly equivalent to the anticanonical class. In fact, we have
$$
p^*(D)=p^*u^*(H_{n-1})=\pi^*(H_{n-1})=h,
$$
which proves the claim. The ramification divisor $R$ can be calculated via the Riemann-Hurwitz formula.
\endproof

\begin{example}{\rm For $n=3$, the quotient $Z:=X_{\lambda}^{4}/<\tau>$ is a del Pezzo surface, as proved in the previous proposition. In particular, the del Pezzo surface $Z$ has degree 2 because $4=h^2=(p^*(-K_Z))^2=2K_Z^2$. The fixed locus of the automorphism $\tau$ is a curve of genus 3 and thus, by the classification of the non--symplectic involution on K3 surfaces, cf. \cite{Nikulin non symplectic}, it follows that $H^2(X_{\lambda}^{4})^{\tau^*}\simeq \langle 2\rangle\oplus\langle -2\rangle^7$, which is indeed the pull-back of the Picard group of $Pic(Z)\simeq \langle 1\rangle\oplus \langle -1\rangle^7$. In particular $X_\lambda^4$ is a subfamily of the $T$-polarized family of K3 surface, for $T\simeq U\oplus U\oplus D_4\oplus \langle -2\rangle^6$. If one considers  a plane sextic curve $C_6$ with 7 ordinary double points, the desingularization of the double cover of $\mathbb{P}^2$ branched along $C_6$ is the general member of the family of $T$-polarized K3 surfaces, and in fact this K3 surface can be constructed as a smooth covering of the blow-up $\widetilde{\mathbb{P}^2}$ of $\mathbb{P}^2$ in the 7 singular points of the sextic. The surface $\widetilde{\mathbb{P}^2}$ is the del Pezzo surface $Z$ and the strict transform of the sextic $C_6$ on $\widetilde{\mathbb{P}^2}$ is the branch locus of the cover $p:X_\lambda^4\ra Z$.}\end{example}

As for the case $n+1$ odd, we can not describe $Z$ in general. We give an example below for $n=4$ in particular we prove
\begin{proposition} The quotient $X_\lambda^5/\tau$ admits a desingularization $Z$ which is a Fano threefold obtained as a double covering of $\mathbb{P}^3$ branched along a reducible sextic.\end{proposition}

Let $\lambda_0 \in {\mathbb C}$ such that $\lambda_0^5 \neq 1$. Consider the Calabi-Yau manifold $X_{\lambda_0}^5$. The restriction $\rho$ to this Calabi-Yau of the projection $\pi$ from $p=(1:-1:0:0:0)$ is generically a degree four covering of ${\mathbb P}^3$, which is not defined at $p$. Let $Bl_p(X_{\lambda_0}^5)$ be the blow-up of $X_{\lambda_0}^5$ at $p$. The map ${\rho}$ lifts to a map $\widetilde{\rho}$ from $Bl_p(X_{\lambda_0}^5)$ to ${\mathbb P}^3$. Let $h$ be the pull-back under the blow-up map of the restriction on $X_{\lambda_0}^5$ of the hyperplane class on ${\mathbb P}^4$. Denote by $e$ the class of the exceptional divisor on $Bl_p(X_{\lambda_0}^5)$. Clearly, the map $\rho$ is associated with the linear system $|h-e|$. Notice that $\rho$ contracts the exceptional divisor $e$.

The automorphism on $X_{\lambda_0}^5$  corresponding to $\tau$ lifts to an automorphism on the blow-up, which fixes $e$ and $h$. The quotient by such an automorphism, say $Z$, is a threefold. By Riemann-Hurwitz's formula, the canonical divisor $K_Z$ satisfies the identity:
\begin{equation}
\label{canonical}
v^*(-K_Y)=h-e,
\end{equation}
where $v: Bl_p(X_{\lambda_0}^5) \rightarrow Z$ is the degree two quotient map. The rational Picard group of $Z$ is generated by $b_1:=u_*(h)$ and $b_2:=u_*(e)$ because $h$ and $e$ belong to a basis of $Pic(Bl_p(X_{\lambda_0}^5))$.

Now, take the divisor on $Z$ given by $\frac{b_1 - b_2}{2}$. By \eqref{canonical} and the push-pull formula, this divisor is equal to $-K_Z$, which is the pull-back of the divisor $h-e$ on the blow-up  $Bl_p(X_{\lambda_0}^5)$. We claim that $h-e$ is ample, so $Z$ is Fano.  To prove that $h-e$ is ample on $Bl_p(X_{\lambda_0}^5)$, let us check the intersection of any irreducible curve and $h-e$. If $C$ comes from a curve that does not pass through $p$, the intersection $(h-e)C$ is positive. If $C$ comes from a curve passing with multiplicity $m$ through $p$, then $(h-e)C$ is positive because the degree of $C$ is greater than $m$. If $C$ lies in the exceptional divisor, then $C$ is linearly equivalent to a multiple of the class $l$ of the line in the exceptional ${\mathbb P}^2$, i.e., $C=dl$ for some positive integer $d$. We have
$$
d(h-e)l=-del=d >0
$$
because $el=-1$. Indeed, $l=c_1\left({\mathcal O}_{{\mathbb P}^2}(1)\right)$ lies in $H^4(Bl_p(X_{\lambda_0}^5))$, so $l=ah^2+be^2$ for some integers $a,b$. Since $hl=0$, we have $l=be^2$. We claim that $b=-1$. As shown in \cite{GH} p. 608, we have
$$
e^2= i_*(i^*(e))=i_*\left(c_1\left({\mathcal O}_{{\mathbb P}^2}(-1)\right)\right)=-l,
$$
where $i$ is the embedding of the exceptional divisor in $Bl_p(X_{\lambda_0}^5)$. Moreover, applying the excess intersection formula and ${\mathcal N}_{l/{\mathbb P}^2}={\mathcal O}_{{\mathbb P}^2}(-1)$, we have
$$
el=i_*(1)l=i_*(i^*(l))=-1;
$$
in particular, $e^3=1$. The degree of $Z$, i.e., the number $(-K_Z)^3$ is $\frac{1}{2}(h^3-e^3)$. Of course, $h^3=5$ because $X_{\lambda_0}^5$ has degree $5$. The morphism associated with the linear system $|-K_Z|$ is a degree two covering of ${\mathbb P}^3$. Indeed, it is easy to check that
$$
\chi(-K_Z)=h^0(Y, -K_Z)=4.
$$

By Riemann-Hurwitz's formula, the ramification divisor is given by $-3K_Z$. If we denote by $\varphi: Z \rightarrow {\mathbb P}^3$ the 2:1 covering above, the branch divisor on ${\mathbb P}^3$ is $nH_3$, where $H_3$ is the hyperplane class on projective three space. Clearly, $\varphi^*(nH_3)=n(-K_Z)=2R=2(-3K_Z)=-6K_Z$, so $n=6$. In other words, the morphism $\varphi$ is branched over a sextic surface in ${\mathbb P}^3$, which is reducible.  We claim that this sextic is the union of a plane and a quintic surface in ${\mathbb P}^3$.

 In fact, let us take into account a generic point on $X_{\lambda_0}^5$ and the point $p$. The line between these two points intersects the Calabi-Yau $X_{\lambda_0}^5$ in three other points. The map $\widetilde{\rho}$ ramifies on the exceptional divisor $e$ and when any of the three points mentioned before coincide. This condition is expressed in terms of the discriminant of a degree three equation. It is a matter of computation to show that the discriminant factors as a product $LC_1C_2^2$, where $L$ is a linear form and $C_i$ is a cubic expression for $i=1,2$ - as homogeneous polynomials on ${\mathbb P}^4$. Hence the map from $X_{\lambda_0}^5$ to ${\mathbb P}^3$ can be described as $w^4=LC_1C_2^2$. Set, now, $z=w^2C_2$; so the covering $Z \rightarrow {\mathbb P}^3$ is given by $z^2=LC_1$. If we introduce homogeneous coordinates $X_0, X_1, X_2, X_3$ in ${\mathbb P}^3$, and use elimination theory, the product $LC_1$ can be written as a reducible sextic that is given by a linear form and a quintic.

\section{The K3 Surfaces $X_{\lambda}^4$}\label{sec4}\label{sec: K3 X}

Little is known about the geometry of $X_{\lambda}^{n+1}$.
Clearly, the case $n=1$ is obvious. The case $n=2$ is the famous
Hesse pencil, \cite{ardo}.
For $n=4$ and $\lambda=1$, the geometry of $X_1^5$ is investigated in \cite{sch}.\\
Here we take into account some aspects of the geometry of $X_{\lambda}^{4}$. In this section, we
denote $X_{\lambda}^4$ by $X_{\lambda}$ for the sake of simplicity.

\subsection{The N\'eron--Severi group of $X^4_{\lambda}$}\label{subsection: The Neron Severi group of X}

Here we want to describe the N\'eron--Severi group of $X_{\lambda}$ considering a special member of the family, whose N\'eron Severi group is known. We need the following results on symplectic automorphisms on K3 surfaces. For each finite group of symplectic automorphisms on a K3 surface (except $G=Q_8$ and $G=T_{24}$) there exists a lattice, $M_G$, which depends only on $G$ and is computed for each $G$ in \cite{Xiao}. As proved in \cite{Nikulin symplectic}, \cite{whitcher} and \cite{Ha}, if the lattice $M_G$ admits a unique primitive embedding in the second cohomology group of a K3 surface, then another lattice, $\Omega_G$, is well defined. It depends only on $G$ and a K3 surface $S$ admits $G$ as group of symplectic automorphisms if and only if $\Omega_G$ is primitively embedded in $NS(S)$. The lattice $\Omega_G$ is $(H^2(S,\Z)^G)^{\perp}$. In the following we will consider the lattice $\Omega_G$ associated to the group $G=\mathfrak{S}_4$, $\mathfrak{A}_4$, $(\Z/4\Z)^2$. For each of these groups there exists a unique embedding of $M_G$ and hence the lattice $\Omega_G$ is well defined.

We consider the surface $X_0$, a special member of the family $X_\lambda$. So we have that $NS(X_0)\supset NS(X_\lambda)$. Moreover we consider the group $H_3 \cong (\Z/4\Z)^2$ acting symplectically on $X_\lambda$ (and hence in particular on $X_0$), and we compute the lattice $\Omega_{H_3}$  as a sublattice of $NS(X_0)$. Then we use the fact that $\Omega_{H_3}\subset NS(X_\lambda)\subset NS(X_0)$ to compute $NS(X_\lambda)$.

Let $F:=X_0$ be the Fermat quartic $\{x_1^4+x_2^4+x_3^2+x_4^4=0\}$. It is a
special member of the family $X_{\lambda}$, the N\'eron--Severi
group $NS(X_{\lambda})\subset NS(F)$. The transcendental lattice of $F$ is known to be $\left[\begin{array}{rr}8&0\\0&8\end{array}\right]$ (cf. \cite{Oguiso}).
On $F$ there are the following 48 lines
$$l_{4a+b}=(\xi_8^{2a+1}s:s:t:\xi_8^{2b+1}t),\ \   l_{4a+b+16}=(s:\xi_8^{2a+1}t:t:\xi_8^{2b+1}s),$$
$$l_{4a+b+32}=(s:t:\xi_8^{2a+1}s:\xi^{2b+1}_8t),$$
where $a=0,1,2,3$, $b=1,2,3,4$ and $\xi_8$ is a primitive $8$-th
root of unity. One can choose a $\Z$-basis of
$NS(F)$ by considering 20 lines among these, indeed one can check that for certain choices of 20 lines among the previous 48 ones, one finds a lattice with discriminant $-8^2$ (so a lattice which is isomorphic to $NS(F)$). In particular, a $\Z$
-basis is given by
$$\{l_1,l_2,l_3,l_4,l_5,l_6,l_7,l_9,l_{10},l_{11},l_{17}, l_{18}, l_{19}, l_{21}, l_{22}, l_{23}, l_{33},l_{34},l_{35}, l_{37} \}.$$

We observe that $l_1\cup l_2\cup l_3\cup
l_4=F\cap\{x_1-\xi_8x_2=0\}$; hence the class $h=l_1+l_2+l_3+l_4$
is the class of the hyperplane section of
$F$.\\

The group $H_3$ acts on the surfaces $X_{\lambda}$ and in particular on the surface $F$. On the surface
$F$ its action transforms the lines generating the N\'eron--Severi group in other lines. Hence, one can completely describe the
action of $H_3$ on the N\'eron--Severi group of $NS(F)$. \\

In particular, the lattice $\Omega_{H_3}(\subset NS(F))$
orthogonal to the invariant lattice $NS(F)^{H_3}$ is generated by the following 18 classes:
$$
\begin{array}{lll}
b_1=-l_{2 }+ l_{37},& b_2=-l_{1 }+ l_{35},& b_3=-l_{2 }+ l_{34},\\
b_4=-l_{1 }+ l_{33},& b_5= -l_{2 }+ l_{23},& b_6=-l_{1 }+ l_{22}, \\
b_7=-l_{2 }+ l_{21},& b_8=-l_{1 }+ l_{19},& b_9=l_{18} - l_{2},
\\
b_{10}= -l_{1 }+ l_{17},& b_{11}=-l_{1 }+ l_{11},& b_{12}=l_{10} -
l_{2}, \\
b_{13}=-l_{1 }+ l_{9},& b_{14}=-l_{2 }+ l_{7},& b_{15}= -l_{1 }+
l_{6},\\
 b_{16}=-l_{2 }+ l_{5},& b_{17}=-l_{2 }+ l_{4},&
b_{18}=-l_{1 }+ l_{3}.
\end{array}
$$

These classes are also contained in $NS(X_{\lambda})$. Indeed, the
symplectic action of $H_3$ on $F$ restricts to a symplectic
action of $H_3$ on $X_{\lambda}$. So $\Omega_{H_3}\subset NS(X_\lambda)$. Moreover, the class $h$
(the class of the hyperplane section) has to be contained in
$NS(X_{\lambda})$ (since $X_{\lambda}$ is the generic member of a
family of quartic hypersurfaces in $\mathbb{P}^3$ on which $H_3$ act symplectically). Hence $\Z h
\bigoplus_{i=1}^{18}\Z b_i\hookrightarrow NS(X_{\lambda})$, where
the inclusion has finite index (because
rank$NS(X_{\lambda})=$rank$(\Z h\oplus \Omega_{H_3})$). Since $NS(X_{\lambda})$ is
primitively embedded in $H^2(X_{\lambda},\Z)\simeq H^2(F,\Z)$,
the lattice $NS(X_{\lambda})$ is primitively embedded in $NS(F)$ and so
$NS(X_{\lambda})\simeq ((\Z h\oplus \bigoplus_{i=1}^{18}\Z
b_i)^{\perp_{NS(F)}})^{\perp_{NS(F)}}.$ A $\Z$-basis of $((\Z
h\oplus \bigoplus_{i=1}^{18}\Z
b_i)^{\perp_{NS(F)}})^{\perp_{NS(F)}}$ in $NS(F)$ is given by the
following 19 effective classes:
$$n_1=h,\ \  \ n_{i}=h+b_{i-1},\
i=2,\ldots, 18,\ \ \  n_{19}=(h+b_{17}+b_{18})/2.$$ We observe
that this shows in particular that $NS(X_{\lambda})$ is an
overlattice of index two of $\Z h\oplus \Omega_{(\Z/4\Z)^2}$,
where
$h$ is the polarization of degree 4 of $X_{\lambda}$.\\

We constructed the lattice $NS(X_{\lambda})$. In particular we can
compute its discriminant group and form:
$NS(X_{\lambda})^{\vee}/NS(X_{\lambda})\simeq (\Z/8\Z)^2\times
\Z/4\Z$ and the discriminant form is
$$M=\left[\begin{array}{lll}\frac{9}{8}&\frac{3}{4}&0\\\frac{3}{4}&\frac{9}{8}&0\\0&0&\frac{5}{4}\end{array}\right].$$
The discriminant form  of the transcendental lattice is the opposite
of the one of the N\'eron--Severi group and its discriminant group
is the same as the one of the N\'eron--Severi group. In particular,
the transcendental lattice $T_{X_{\lambda}}$ is a rank
3(=$22-\rho(X_{\lambda})$) lattice with signature $(2,1)$,
discriminant group $(\Z/8\Z)^2\times \Z/4\Z$ and with discriminant
form $-M$. We recall the following trivial proposition.

\begin{proposition}\label{prop: condition on T=L(2)}
Let $(T,b_T$) be a lattice such that: {\it i)} rk$T=l(T)$;\\ {\it
ii)} $T^{\vee}/T\simeq \bigoplus_i \Z/2d_i\Z$ and $\beta_i$
generate $\Z/2d_i\Z$ in $T^{\vee}/T$.\\
Let $L$ be the free $\Z$-module $L=\{x\in T\otimes \Q\mbox{ such
that }2x\in T\}$ and $b_L=2(b_{T\otimes\Q})_{|L}$ a bilinear form
defined on $L$.
Then $(L,b_L)$ is a lattice and $T\simeq L(2)$ if and only if $b_{T^\otimes \Q}(d_i\beta_i,d_j\beta_j)\in \frac{1}{2}\Z$.\\
Moreover, $L$ is an even lattice if and only if
$b_{T\otimes\Q}(d_i\beta_i,d_i\beta_i)\in \Z$.
\end{proposition}

Let us suppose that $T\simeq L(2)$. Then rk$(L)$=rk$(T)$,
$d(L)=d(T)/2^{rk(T)}$. With the same notations of Proposition
\ref{prop: condition on T=L(2)}, the discriminant group of $L$ is
generated by $\lambda_i=2\beta_i$, for $i=1,\ldots n$ such that
$2\beta_i\notin L$ and its discriminant form is 2 times the
discriminant form of $T$.
\begin{proposition} Let $S$ be a K3 surface such that $T_S\simeq L(2)$ for an even lattice $L$. If $L$ admits a primitive embedding in $U\oplus U\oplus U$, then $S$ is a Kummer surface.\end{proposition}
\proof Let $A$ be the Abelian surface such that $T_A\simeq L$ (this surface exists by the condition on $L$). The transcendental lattice of $Km(A)$ is $T_{Km(A)}\simeq T_A(2)\simeq L(2)$ (cf. \cite{morrison}). In particular $T_{Km(A)}\simeq T_S$. Since $L\subset U\oplus U\oplus U$, $rk(L)\leq 6$. The N\'eron--Severi group of a K3 surface with transcendental lattice of rank less then or equal to 6, is uniquely determined by the transcendental lattice. Hence $NS(S)\simeq NS(Km(A))$. Since $Km(A)$ is a Kummer surface there exists a primitive embedding of the Kummer lattice $K$ in $NS(Km(A))$, and so there is a primitive embedding of $K$ in $NS(S)$. This implies that $S$ is a Kummer surface (cf. \cite{nikulin kummer}).\endproof

\begin{corollary} The surface $X_{\lambda}$ is a Kummer surface for any $\lambda$.
\label{kummer}
\end{corollary}

\proof For $\lambda$ such that $\lambda^4=1$, this is known: see,
for instance, \cite{hud}. As for the other values of $\lambda$, it
is easy to verify that $T_{X_{\lambda}}$ satisfies the condition
of Proposition \ref{prop: condition on T=L(2)}; hence there exists
a lattice $R$ such that $T_{X_{\lambda}}\simeq R(2)$. The lattice
$R$ is even, of rank 3, signature $(2,1)$ and its discriminant
form is given by $-2M$.

We recall that each even lattice of signature $(2,1)$ admits a
primitive embedding in $U\oplus U\oplus U$, and hence there exists
an Abelian surface $A$ such that $T_A\simeq R$. This is enough to
prove that $X_{\lambda}$ is the Kummer surface of an Abelian
surface $A$.\endproof

\subsubsection{Automorphisms of the K3 Surfaces $X_{\lambda}^4$}
We already observed (Prop. \ref{prop: An and Sn groups of autom}) that the maximal finite group  $G_{L}$ of symplectic
automorphisms of $X_\lambda$, which come from projective transformations of ${\mathbb P}^3$ is isomorphic to $\mathfrak{A}_4\rtimes H_3$ for a generic $\lambda$.
\begin{remark} Let $G$ be the maximal finite group of symplectic
automorphisms of the K3 surface $X_{\lambda}$. If $\lambda=-3$,
then $G=(\Z/4\Z)^2\rtimes \mathfrak{A}_5$ (cf. \cite{Mu}). If $\lambda=0$, then $G=F_{384}$ (cf. \cite{Oguiso}).\\
Let $G_L$ be the maximal finite group of symplectic automorphisms
of $X_{\lambda}$ preserving the polarization $L$ realizing
$X_{\lambda}$ as quartic hypersurface in $\mathbb{P}^3$ with
equation $\sum_{i=1}^4x_i^4-4\lambda \prod_{i=1}^4x_i=0$. If
$G_L=(\Z/4\Z)^2\rtimes {\mathfrak A}_5$, then $\lambda=-3$ (cf. \cite{Mu}). If $G_L=F_{384}$, then $\lambda=0$ (cf.
\cite{Oguiso}).
\end{remark}

Since for a generic $\lambda$ $\mathfrak{A}_4$ is a group of symplectic automorphism of $X_\lambda$ and $\mathfrak{S}_4$ is not, $\Omega_{\mathfrak{A}_4}\subset NS(X_{\lambda})$ and $\Omega_{\mathfrak{S}_4}\not\subset NS(X_\lambda)$. In our context, this follows directly also from the construction of $NS(X_\lambda)$ as sublattice of $NS(F)$ we considered. Indeed we chose a
class $v\in NS(F)$ such that ${\mathbb Z}v:=(\Z h \bigoplus_{i=1}^{18} \Z
b_i)^{\perp_{NS(F)}}$ and we said that
$NS(X_{\lambda})=v^{\perp_{NS(F)}}$. By explicit computation one
sees that the class $v \equiv -l_2-l_4+l_9+l_{11}$ in $NS(F)$. Notice that $$v\notin
NS(F)^{\mathfrak{S}_4},\ \mbox{hence} \ \Omega_{\mathfrak{S}_4}=
(NS(F)^{\mathfrak{S}_4})^{\perp_{NS(F)}}\not\subset
v^{\perp_{NS(F)}}=NS(X_{\lambda})$$ but $$v\in
NS(F)^{\mathfrak{A}_4} \mbox{ hence }\Omega_{\mathfrak{A}_4}=
(NS(F)^{\mathfrak{A}_4})^{\perp_{NS(F)}}\subset
v^{\perp_{NS(F)}}=NS(X_{\lambda}),$$ so we obtain
$\Omega_{\mathfrak{A}_4}\subset NS(X_{\lambda})$ and
$\Omega_{\mathfrak{S}_4}\not\subset
NS(X_{\lambda})$.

In the Section \ref{subsection: The Neron Severi group of X} we gave an explicit description of the
generators of $NS(F)$. In Remark \ref{rem: autom X_0}, we describe the symplectic action of $\mathfrak{S}_4$ on $F=X_0$. Hence we can explicitly compute
$\Omega_{\mathfrak{S}_4}\simeq H^2(F,\Z)^{\mathfrak{S_4}}$,
recalling that the action of $\mathfrak{S}_4$ is trivial on the
transcendental lattice of $\mathfrak{S}_4$. We have that
$\Omega_{\mathfrak{S}_4}$ is generated by
$$
\begin{array}{llll}
d_1=l_{37} - l_{5},& d_{2}= l_{2} - l_{22} - l_{23} + l_{35}, \\
d_{5}=
 -l_{1} + l_{17},&  d_{6}=
 -l_{17} - l_{2} + l_{22} + l_{5} \\
d_{3}= -l_{2} + l_{34},& d_{4}= -l_{1} + l_{33},\\
d_{7}= -l_{2} + l_{21} + l_{23} - l_{4},& d_{11}= -l_{1} + l_{11},\\
 d_{9}= l_{18} - l_{2},&d_{10}= l_{1} -
l_{17} - l_{21} + l_{5}, \\
d_{8}= l_{19} - l_{21} - l_{35} +
l_{37},& d_{12}= l_{10} -
 l_{4}, \\
 d_{13}= -l_{35} + l_{9},& d_{14}= -l_{2} - l_{4} + l_{5} + l_{7}, \\ d_{15}=
 -l_{1} - l_{34} + l_{37} + l_{6},& d_{16}= l_{3} - l_{35}\\
 d_{17}=-l_{1} + l_{37}
\end{array}
$$ and $\Omega_{\mathfrak{A}_4}$ is generated by the classes
$d_i$, $i=1,\ldots,16$.\\
In particular $\Omega_{\mathfrak{S}_4}$ is a rank 17 negative
definite lattice with discriminant group $\Z/4\Z\times
(\Z/12\Z)^2$ and $\Omega_{\mathfrak{A}_4}$ is a rank 16 negative
definite lattice with discriminant group $(\Z/2\Z)^2\times
(\Z/12\Z)^2$. Both these lattices are generated by vectors of
length $-4$, which is their maximal length.\\

In \cite{dihedral}, it is observed that there exists pair $(G,H)$ of groups acting symplectically on K3 surfaces such that $H$ is a subgroup of $G$ and a K3 surface admits $G$ as group of symplectic automorphisms if and only if it admits $H$ as group of symplectic automorphisms. If this property holds for the pair $(G,H)$ and the lattices $\Omega_G$ and $\Omega_H$ are both well defined, then $\Omega_G\simeq \Omega_H$. A complete list of the pair of groups $(G,H)$ with this property can be found in \cite{Ha}, but in general no explicit examples are known. The pair $(\mathfrak{A}_4\rtimes (\Z/4\Z)^2, (\Z/4\Z)^2)$ has this property and the K3 surface $X_\lambda^4$ is an explicit geometric example. In particular the lattice $\Omega_{\mathfrak{A}_4\rtimes (\Z/4\Z)^2}$ is isometric to $\Omega_{(\Z/4\Z)^2}$, which is computed in \cite{symplectic non prime}: it is a negative defined lattice of rank 18, with discriminant group $(\Z/2\Z)^2\times (\Z/8\Z)^2$ and it is generated by vectors of maximal length, i.e. of length $-4$.\\

We observed that a K3 surface admits $(\Z/4\Z)^2$ as a group of symplectic automorphisms, then it
admits also $\mathfrak{A}_4$ as a group of symplectic automorphisms, in fact admits $\mathfrak{A}_4\rtimes (\Z/4\Z)^2$.
If we restrict our attention only to hypersurfaces of degree $n+1$
in $\mathbb{P}^n$ and we fix a particular action of
$(\Z/(n+1)\Z)^{n-1}$ (the one given by $H_n$), this is true also
for higher dimension, indeed (by Propositions \ref{prop:
Z/(n+1)Z^(n-1) group of autom}, \ref{prop: An and Sn groups of autom}) the family of hypersurfaces of degree $n+1$ in
$\mathbb{P}^n$ admitting $H_n$ as a group of symplectic
automorphisms admits also $\mathfrak{A}_{n+1}$ as group of
symplectic automorphisms. It would be interesting to construct a
Calabi-Yau manifold $V$ of dimension $n-1>2$ such that $V$ admits
$H_n$ but not $\mathfrak{A}_{n+1}$ as a group of automorphisms preserving the period.

\section{The Quotients of the CY Three-folds $X_{\lambda}^5$}\label{sec5}\label{sec: Quotients CY 3 fold}

Let $V$ be a Calabi--Yau variety of dimension at most 3 and let $G$ be a finite group of automorphisms of $V$ which preserves the period of $V$. If the quotient $V/G$ is smooth, then it is a Calabi--Yau variety. If it is singular, then there exists a desingularization which is a Calabi--Yau variety, cf.\ \cite{yas}. Such a desingularization is not unique if the dimension of $V$ is 3, but the Hodge numbers of all the desingularizations of $V/G$ which are Calabi--Yau varieties are the same \cite{Bat2}. These numbers can be computed via orbifold cohomology (\cite{CheR}). The aim of this section is to consider some quotients of $X_\lambda^5$ by groups of automorphisms preserving the period. We recall that the full group of automorphism of $X_\lambda^5$ is $ \mathfrak{S}_{5}\rtimes (\Z/5\Z)^{3}$ and the subgroup which acts trivially on the period is $\mathfrak{A}_{5}\rtimes (\Z/5\Z)^{3}$.
\subsection{Automorphisms Fixing the Period of Calabi--Yau Three-folds}
First, we recall some known relations among the fixed locus of a group of automorphisms and the action of such a group on the cohomology. Let $J$ be a group of automorphisms on a variety $Z$, and $j\in J$, then:
\begin{equation}\label{formula: Lefschetz fixed points}\sum_{k=0}^{2dim(Z)}(-1)^k tr(j_{|H^k(Z,\Q)})=\chi(Z^h), \mbox{ (Lefschetz fixed points formula)}\end{equation}
and by \cite{FH}, p. 16,
\begin{equation}\label{formula: rk trivial on cohom}rk(H^*(Z,\Z)^J)=\frac{1}{|J|}\sum_{j\in J}tr(j_{H^*(X,\Q)})=\frac{1}{|J|}\sum_{j\in J}\chi(X^j).\end{equation}
From now on let $G$ be a group of automorphisms fixing the period of a Calabi--Yau threefold $Z$.

Let us denote by $p_{1,1}:=dim(H^{1,1}(Z)^G)$ and $p_{1,2}:=dim(H^{1,2}(Z)^G)$. Then the formula \eqref{formula: rk trivial on cohom} becomes
$$
2(1+p_{1,1})-(2+2p_{1,2})=\frac{1}{|G|}\sum_{g\in G}\chi(X^g).
$$

We observe that, by definition, $G$ acts trivially on $H^{3,0}$ but its action on $H^{2,1}$ and on $H^{1,1}$ depends on the group and on the threefold. Howerver, under certain conditions on the group of automorphisms or on the threefold, one can assume that the action either on $H^{1,1}$ or on $H^{2,1}$ is trivial.\\

If $h^{1,1}(Z)=1$, we have already noticed in the proof of Proposition \ref{prop: An and Sn groups of autom}, that the action of each automorphism on $H^{1,1}(Z)$ is trivial and thus Formula \eqref{formula: rk trivial on cohom} allows us to deduce $p_{1,2}$ by the Euler characterisitc of the fixed locus of the elements in $G$:
\begin{equation}\label{eq: p12}p_{1,2}=1-\frac{1}{2|G|}\sum_{g\in G}\chi(X^g).\end{equation}
In particular, \eqref{eq: p12} applies to every subgroups of $\mathfrak{A}_5\rtimes H_4$ acting on $X_\lambda^5$.\\

An automorphism on a Calabi--Yau threefolds is called {\it maximal} if it extends to an automorphism of the family of deformations of the Calabi--Yau and restricts to an automorphism of each fiber. A maximal automorphism of a Calabi--Yau three-fold which preserves the period acts trivially on all of the middle cohomology, in particular on $H^{2,1}$ and $H^{1,2}$.
\begin{example}{\rm The Calabi--Yau three-fold $X_\lambda^5$ does not admit maximal automorphisms, indeed the generic hypersurface of degree 5 in $\mathbb{P}^4$ is a member of the deformation family of $X_\lambda^5$ and does not admit automorphisms.\\
The group $\mathfrak{A}_5$ is a maximal group of automorphism which preserves the period of $Y_\lambda^5$. Indeed $h^{2,1}(Y_\lambda^5)=1$, thus all the deformations of $Y_\lambda^5$ depend  on the variation of the parameter $\lambda$, i.e. all the member of the family have an equation of type $(\sum_{i=1}^5y_i)^5+t\prod_{i=1}^5y_i$ for a certain $t$, and hence admit $\mathfrak{A}_5$ as automorphism group fixing the period. }\end{example}

In order to compute the Hodge numbers of a crepant resolution of a quotient $X_\lambda^5/G$ for a certain subgroup $G$ of $\mathfrak{A}_5\rtimes H_4$, one can apply the results in \cite{CheR}, where the orbifold cohomology is defined and it is proved that, under some conditions, it coincides with the ordinary cohomology of a desingularization. To this end, we recall how to compute the orbifold cohomology of the quotient of a Calabi--Yau threefold $Z$ by a group $G\subset Aut(Z)$ of automorphisms that fix the period of $Z$.

Let $S$ be the set of representatives of the conjugacy classes of $G$. For each $s\in S$, let $F_s$ be the fixed locus $Fix_s(Z)=\{z\in Z\ |\ s(z)=z\}$; denote by  $F^i_s$ the connected components of $F_s$. For each component $F^i_s$ the age of $s$ is well defined and does not depend on the representative in $S$. Let $C_s$ be the centralizer of $s$ in $G$. Then the orbifold cohomology of $Z/G$ is
\begin{equation}\label{eq:orbifold cohom}H^{p,q}(Z/G)_{orb}=H^{p,q}(Z)^G\bigoplus \bigoplus_{s\in S,\ s\neq 1}\bigoplus_i H^{p-age(s),q-age(s)}(F_s^i/C_s).\end{equation}
In order to apply the orbifold cohomology, one has to compute the age of the elements in $g\in G$ near the fixed locus. The age of an automorphism acting on a variety near the fixed locus depends on the diagonalization of the action of the automorphism, and so in general both on the automorphism and on the variety. Anyway, if one considers an automorphism $g$ fixing the period of a Calabi--Yau threefold $Z$, the Hodge numbers of a crepant resolution of $Z/g$ depend only on the topology of the fixed locus and on the dimension of the spaces $H^{1,1}(V)^g$ and $H^{1,2}(V)^g$. In particular, the following holds:
\begin{proposition}\label{cor: Hodge numbers quotient} Let $p$ be a prime. Let $Z$  be a Calabi--Yau threefold with $h^{1,1}(Z)=1$, $g$ be an automorphism of order $p$ of $Z$ preserving the period and $F_g=\coprod_{i=0}^m P_i\coprod_{j=0}^k C_j.$ Let $\widetilde{Z/g}$ be a crepant resolution of the quotient $Z/g$. Then $\widetilde{Z/g}$ is a Calabi--Yau variety and the following holds: $$h^{1,1}(\widetilde{Z/g})=1+\frac{p-1}{2}m+(p-1)k,$$
$$h^{1,2}(\widetilde{Z/g})=1-\frac{1}{2p}\left(\chi(Z)+(p-1)\left(m+2k-2\sum_{i=1}^kg(C_i)\right)\right)+(p-1)\sum_{i=1}^kg(C_i).$$
\end{proposition}
\proof In order to compute the Hodge numbers of a crepant resolution of $Z/g$ we use the formula \eqref{eq:orbifold cohom}. Thus, we have to consider the age of $g^i$ near the fixed locus for $i=1,\ldots p-1$. The fixed locus of $g^i$ coincides with that of $g$ because the order of $g$ is a prime number. If $F$ is a connected component of the fixed locus of $g$, the age of $g^i$ near $F$ satisfies the following properties (cf. \cite[Lemma 3.2.1]{CheR})
 \begin{equation}
 age(g^i)+age(g^{-i})=codim(F), \ \ age(g^j)\in\N_{>0},\ \ j=1,\ldots (p-1).
 \label{agegginverse}
 \end{equation}
Let $F=P_i$ be an isolated fixed point, then $age(g^i)+age(g^{-i})=3$ and thus $\{age(g^i),age(g^{-i})\}=\{1,2\}$. We recall that $g^{-i}=g^{p-i}$ and thus the set $\{g^i\}_{i=1,\ldots p-1}$ coincides with the set $\{g^i,g^{-i}\}_{i=1,\ldots (p-1)/2}$. Thus for every isolated fixed point $P\in Fix_g(Z)$ there exists exactly $(p-1)/2$ elements in $\langle g \rangle$ such that their age near $P$ is 1 and exactly $(p-1)/2$ elements in $\langle g \rangle$ such that their age near $P$ is 2. Similarly, if $F=C$ is a fixed curve, every element in $\langle g \rangle$ has age $1$ near $C$. The statement follows from \eqref{eq: p12} and \eqref{eq:orbifold cohom}.\endproof

In some particular cases (e.g., $|g|=2,3$), one can describe more directly the desingularization of the quotient of $Z/g$. Indeed, the type of singularities of a quotient of a Calabi--Yau three-fold by an automorphism is a local property and its desingularization depends only on the local action of $g$ on $Z$ near the fixed locus.\\ If $|g|=2$, the local action of $g$ is given by the diagonal matrix $diag(-1,-1,1)$. In this case it suffices to blow-up the fixed curve in $Z$, introducing a $\mathbb{P}^1$-bundle over the fixed curve and obtaining the variety $\widetilde{Z}$. The involution $g$ lifts to $\widetilde{g}$ on $\widetilde{Z}$ and a local computation shows that $\widetilde{g}$ fixes the exceptional divisor over the fixed curve. Thus,  the quotient $\widetilde{Z}/\widetilde{g}$ is a smooth threefold, which is a desingularization of $Z/g$. Since the fixed curves are disjoint, one can blow-up each of them independently. Let us denote by $E_i$, $i=1,\ldots k$ the exceptional divisors of the simultaneous blow-up of the fixed curves; let $\beta:\widetilde{Z}\ra Z$ be the blow-up, $\pi:\widetilde{Z}\ra \widetilde{Z}/\widetilde{g}$ be the quotient map. The variety $\widetilde{Z}/\widetilde{g}$ is a Calabi--Yau variety: indeed,  $K_{\widetilde{Z}}=\beta^*(K_Z)+E$ and $K_{\widetilde{Z}}=\pi^*K_{\widetilde{Z}/\widetilde{g}}+E$, where $E$ is the contribution given by the exceptional divisors; thus $\pi^*K_{\widetilde{Z}/\widetilde{g}}$ is trivial.

\begin{remark} Since the construction of the desingularization of the quotient $Z/g$ depends only on the local action of $g$ near the fixed locus and on the fact that the canonical bundle of $Z$ is trivial, one can deduce the same result observing that our construction is the same applied by \cite{V2} to obtain a desingularization of $(S\times E)/\iota_S\times \iota_E$ where $S$ is a K3 surface, $E$ is an elliptic curve (in fact the canonical bundle of $S\times E$ is trivial) and $\iota_S\times\iota_E$ is an involution which fixes curves and  near the fixed locus acts as $diag(-1,-1,1)$.
\end{remark}
Similarly, if $|g|=3$, one can construct the desingularization of  $Z/g$ by blowing up $Z$ and take the quotient of the blow-up by the automorphism induced by $g$. By the computation of the Hodge numbers, it is clear that we introduce one exceptional divisor for each point and 2 for each curve in the fixed locus. Indeed it suffices to blow up $Z$ once in the fixed points and to blow up three times for every fixed curves. In the case of the fixed curves one has to contract one of the exceptional curves in order to obtain a minimal smooth Calabi--Yau three-fold, thus yielding two exceptional divisor over a fixed curve. The explicit computation is shown in \cite{Rohde}, where the desingularization of a quotient of the 3-fold $S\times E_{\xi_3}$ is constructed, where $S$ is a K3 surface and $E_{\xi_3}$is an elliptic curve, and the quotient is taken w.r.t an order three group preserving the period.

\subsubsection{Quotients of $X_\lambda^5$ by the Cyclic Groups $\Z/2\Z$, $\Z/3\Z$, $\Z/5\Z$}\label{sec: quotient cyclic groups}
Here we compute the Hodge numbers of some quotients of $X_\lambda^5$ applying the results of the previous section. Thus we will apply Proposition \ref{cor: Hodge numbers quotient} to the case $Z=X_\lambda^5$ and to certain automorphisms $g\in \mathfrak{A}_5\rtimes H_4$.\\
Let us consider the particular case $g=(12)(34)$: so $|g|=2$ and the fixed locus of $g$ on $X_\lambda^5$ consists of a line, $(s:-s:t:-t:0)$, and of a plane smooth quintic, $2x_1^5+2x_3^5+x_5^5+5\lambda x_1^2x_3^2x_5=0$. Thus each desingularization of $X_\lambda^5$ which is a Calabi--Yau 3-fold has the following Hodge numbers
$$h^{0,0}=h^{3,0}=h^{0,3}=1,\ h^{1,0}=h^{0,1}=h^{2,0}=h^{0,2}=0,\ h^{1,1}=3,\ h^{2,1}=59.$$
\begin{remark} The dimension of the eigenspace $H^{2,1}(X_\lambda^5)_1$ is 53. This can be also deduced recalling that $H^{2,1}$ parametrizes the deformations of $X_\lambda^5$, which is a quintic in $\mathbb{P}^4$. Thus, it can be identified with the space of the monomials of degree 5 in 5 variables modulo the projective transformations of $\mathbb{P}^4$. Hence, the eigenspace relative to the eigenvalue 1 for the automorphisms $(12)(34)$ can be identified with the space of the monomials among the previous ones which are invariant under $(12)(34)$.
\end{remark}

Similarly one computes the following Hodge numbers of a crepant resolution $\widetilde{X_\lambda^5/g}$ of the quotient of $X_\lambda^5$ by certain automorphisms $g\in\mathfrak{A}_5\rtimes H_4$. Since $\widetilde{X_\lambda^5/g}$ is a Calabi--Yau threefold, the only numbers $h^{1,1}$ and $h^{2,1}$ are given (the others do not depend on $g$):

$$\begin{array}{|c|c|c|c|}
\hline
g&|\langle g\rangle|&h^{1,1}(\widetilde{X_\lambda^5/g})&h^{2,1}(\widetilde{X_\lambda^5/g})\\
\hline
g=(12)(34)&2&3&59\\
\hline
g=(123)&3&5&49\\
\hline
g=h_{(0,1,4,0,0)}&5&5&49\\
\hline
g=h_{(0,1,1,3,0)}&5&21&17\\
\hline
g=h_{(0,1,2,3,4)}&5&1&21\\
\hline
g=(12345)&5&1&21\\
\hline
\end{array}$$

\begin{remark} Since both $h_{(0,1,2,3,4)}$ and $(12345)$ are automorphisms of order 5 fixing the period and fixed points free, the Hodge numbers of the smooth three-folds $X_\lambda^5/h_{(0,1,2,3,4)}$ and $X_\lambda^5/(12345)$ are clearly the same. In both these cases the fundamental group is $\Z/5\Z$.\\
It is more surprising that the Hodge numbers of a crepant resolution of $X_\lambda^5/(123)$ coincides with the ones of a crepant resolution of $X_\lambda^5/h_{(0,1,4,0,0)}$. It would be interesting to understand whether there is a deeper and geometric reason for this phenomenon.
\end{remark}

\subsection{Quotients by Subgroups of $\mathfrak{A}_5$}\label{sec: quotient subgroups A5}

In \cite{St}, quotients of the Fermat quintic three-fold $F_3^5\subset \mathbb{P}^4$ are considered. The automorphisms of $\mathfrak{A}_5$ act both on $F_3^5=X_0^5$ and on $X_\lambda^5$. Moreover, for each automorphism in $\mathfrak{A}_5$ the fixed locus of its action on $F_5^3$ and on $X_\lambda^5$ has the same number of fixed points, the same number of fixed curves and the fixed curves have the same genus. Similarly for each automorphism in $\mathfrak{A}_5$, the diagonalizations of its action  on $F_3^5$ and on $X_\lambda^5$ are the same. Since the Hodge numbers of a crepant resolution of the quotient of a Calabi--Yau threefolds by a group $G$ depend only on the properties of the fixed locus and of the age of every automorphisms in $G$,  by \cite{St} we obtain that:
$$\begin{array}{|l|c|c|}
\hline
Group&h^{11}(\widetilde{X_\lambda^5/G})&h^{2,1}(\widetilde{X_\lambda^5/G})\\
\hline
\mathfrak{A}_5&5&15\\
\hline
\mathfrak{A}_4=\langle(12)(34),(123)\rangle&7&29\\
\hline
\mathcal{D}_5=\langle(12)(35),(12345)\rangle&3&19\\
\hline
\mathfrak{S}_3=\langle(12)(45),(23)(45)\rangle&5&33\\
\hline
\Z/5\Z=\langle(12345)\rangle&1&21\\
\hline
(\Z/2\Z)^2=\langle (12)(34),(13)(24)\rangle&7&41\\
\hline
\Z/3\Z=\langle(123)\rangle&5&49\\
\hline
\Z/2\Z=\langle(12)(34)\rangle&3&59\\
\hline
\end{array}$$
The cases where $G$ is cyclic were already considered in the previous section.

\subsection{Quotients by Subgroups of $H_4$}

The quotient of $X_\lambda^5$ by $H_4=(\Z/5\Z)^3$ is well known, and it is described in Section \ref{sec: Y5}.
We recall that the automorphisms  in $H_4$ are of the form $h_{(a_1,a_2,a_3,a_4,a_5)}:(x_1:x_2:x_3:x_4:x_5)\ra (\xi^{a_1} x_1:\xi^{a_2} x_2:\xi^{a_3} x_3:\xi^{a_4} x_4: \xi^{a_5}x_5)$ and we can always assume $a_i=0$ for at least one $i\in\{1,2,3,4,5\}$, since $h_{(a_1,a_2,a_3,a_4,a_5)}$ is an automorphism of a projective space. From the point of view of the fixed locus, an automorphism of order 5 in $(\Z/5\Z)^3$ is of one of the following type:\\
$i)$ the automorphism $h_{(0,a_2,a_3,a_4,a_5)}$ fixes a smooth plane curve of degree 5 (thus of genus 6) and is induced on $X_\lambda^5$ by an automorphism with exactly two values in $\{a_i\}$ $i=2,3,4,5$  which are not trivial mod $5$;\\
$ii)$ the automorphism $h_{(0,a_2,a_3,a_4,a_5)}$ fixes 10 points and is induced on $X_\lambda^5$ by an automorphism with exactly three values in $\{a_i\}$ $i=2,3,4,5$  which are not trivial mod $5$ (this implies that two of these three values are the equal);\\
$iii)$ the automorphism $h_{(0,a_2,a_3,a_4,a_5)}$ is fixed points free and  is induced on $X_\lambda^5$ by an automorphism with $a_i\not\equiv 0\mod 5$ $i=2,3,4,5$.\\

We already computed the Hodge numbers of a crepant resolution, $\widetilde{X_\lambda^5/g}$, of the quotients of $X_\lambda^5$ by cyclic subgroups of $H_4$ in Section \ref{sec: quotient cyclic groups}.\\

The subgroups $(\Z/5\Z)^2\subset (\Z/5\Z)^3$ are of three different types with respect to the set of points with a non trivial stabilizer: \\
i) the group $G_1:=\langle h_{(1,4,0,0,0)},h_{(0,0,1,4,0)}\rangle$, \\
ii) the group $G_2:=\langle h_{(1,4,0,0,0)},h_{(1,0,4,0,0)}\rangle$, \\
iii) and the group $G_3:=\langle h_{(1,1,0,0,3)},h_{(1,3,1,0,0)}\rangle$.\\

All the other subgroups $G\simeq (\Z/5\Z)^2$ are of these types, in the sense that they contain the same numbers of elements with the same type of fixed locus, and thus the crepant resolutions of the quotient $X_\lambda^5/G$ have the Hodge numbers of the crepant resolutions of one of the quotients $X_\lambda^5/G_i$ for $i=1,2,3$.

Let us compute the Hodge numbers of a crepant resolution of $X_\lambda^5/G_1$; the other computations are very similar.
The elements of the group $G_1$ are
$$id,  h_{(i,4i,0,0,0)}, h_{(0,0,i,4i,0)}, h_{(i,4i,i,4i,0)}, h_{(i,4i,2i,3i,0)}, h_{(i,4i,3i,2i,0)}, h_{(i,4i,4i,i,0)}$$ for $i=1,2,3,4$. The 8 elements $h_{(i,4i,0,0,0)}$, $h_{(0,0,i,4i,0)}$, $i=1,2,3,4$ fix a curve of genus 6 and their age near the fixed locus is 1. The 8 elements $h_{(i,4i,i,4i,0)}$, $h_{(i,4i,4i,i,0)}$, $i=1,2,3,4$ fix 10 points and 4 of them have age 1 near the fixed locus, 4 have age 2. The 8 elements $h_{(i,4i,2i,3i,0)}$, $h_{(i,4i,3i,2i,0)}$ act freely.
Thus $rk(H^3(X_\lambda^5)^G_{1})=12$, by \eqref{formula: rk trivial on cohom}.

In order to compute the orbifold cohomology of the quotient $X_\lambda^5/G_1$, we need to compute the quotients of $Fix_g(X_\lambda^5)/G_1$ for every $g\in G_1$. For example, let $C$ be the curve fixed by $h_{(1,4,0,0,0)}$: the element $h_{(0,0,1,4,0)}$ is an automorphism of $C$ without fixed points, thus $C/G_1\simeq C/h_{(0,0,1,4,0)}$ is a curve of genus 2. Similarly one proves that $Fix_{h_{(0,0,i,4i,0)}}(X^5_\lambda)/G_1$ is a curve of genus 2, $Fix_{h_{(i,4i,i,4i,0)}}(X_\lambda)^5/G_1$ and $Fix_{h_{(i,4i,4i,i,0)}}(X^5_\lambda)/G_1$ consist of 2 points. Thus one finds:
$$h^2(\widetilde{X^5_\lambda}/G_1)=1+8+4(2)=17,$$ and $h^3(\widetilde{X_\lambda^5}/G_1)=12+8(4)=44$. In particular, $$h^{1,1}(\widetilde{X_\lambda^5}/G_1)=17, h^{2,1}(\widetilde{X_\lambda^5}/G_1)=21.$$
Similarly, one computes the other Hodge numbers and finds:
$$\begin{array}{|c|c|c|}
\hline
G_1=\langle h_{(1,4,0,0,0)},h_{(0,0,1,4,0)}\rangle&h^{1,1}(\widetilde{X_\lambda^5/G_1})=17& h^{2,1}(\widetilde{X_\lambda^5/G_1})=21\\
\hline
G_2=\langle h_{(1,4,0,0,0)},h_{(1,0,0,4,0)}\rangle&h^{1,1}(\widetilde{X_\lambda^5/G_2})=49& h^{2,1}(\widetilde{X_\lambda^5/G_2})=5\\
\hline
G_3=\langle h_{(1,1,0,0,3)},h_{(1,3,1,0,0)}\rangle&h^{1,1}(\widetilde{X_\lambda^5/G_3})=21& h^{2,1}(\widetilde{X_\lambda^5/G_3})=1\\
\hline
\end{array}$$

It is easy to show that there are no groups $G\simeq (\Z/5\Z)^2$ such that a crepant resolution of $X_\lambda^5/G$ has Hodge numbers different from those of $X_\lambda^5/G_j$ for a certain $j=1,2,3$. In fact, the generators of each $(\Z/5\Z)^2$ have order 5 in $(\Z/5\Z)^3$ and there are only three types of such automorphisms, with respect to the fixed locus. One has to chose two automorphisms among them as generators of the group $G$. Let us choose $h_{(1,4,0,0,0)}$ as the first generator: the second one could be of the same type, but with 1 and 4 in a different position (this gives the groups $G_1$ and $G_2$) or of type $h_{(a,b,c,d,e)}$ where $(a,b,c,d,e)$ is a permutation of $(1,1,0,0,3)$  or of type $h_{(a,b,c,d,e)}$ $\{a,b,c,d,e\}=\{0,1,2,3,4\}$. We observe that the last case gives the group $G_1$, since in $G_1$ there exists an element $h_{(a,b,c,d,e)}$ with $\{a,b,c,d,e\}=\{0,1,2,3,4\}$. In case the second generator is of type $h_{(a,b,c,d,e)}$, where $(a,b,c,d,e)$ is a permutation of $(1,1,0,0,3)$, a longer but similar analysis shows that one obtains again either $G_1$ or $G_2$. Similarly, one finds that the unique other possibility is the group $G_3$, if the group does not contain an element $h_{(a,b,c,d,e)}$, where $(a,b,c,d,e)$ is a permutation of $(1,4,0,0,0)$.
\begin{remark} The Hodge numbers of a crepant resolution of $X_\lambda^5/G_1$ (resp.  $X_\lambda^5/G_2$,  $X_\lambda^5/G_3$) are mirror of the ones of a crepant resolution of $X_\lambda^5/h_{(0,1,1,3,0)}$ (resp.  $X_\lambda^5/h_{(0,1,4,0,0)}$, $X_\lambda^5/h_{(0,1,2,3,4)}$). This is an example of the correspondence of \cite[Theorem 4]{ChiR} putting $G:=G_1$ and $G^T:=\langle h_{(0,1,1,3,0)}\rangle$ (cf. \cite[Equation (6) and Theorem 4]{ChiR}).\\
Let $X_1$ be the quotient $(X_\lambda^5)/h_{(1,4,1,4,0)}$ and $h^1_{(a,b,c,d,e)}$ the automorphism induced by $h_{(a,b,c,d,e)}$ on $X_1$. There is a chain of degree 5 quotients $X_\lambda^5\ra X_1\ra X_2:=X_1/h^1_{(1,4,0,0,0)}\ra Y_\lambda^5$ such that $X_2\simeq X_\lambda^5/G_1$ and the crepant resolutions of $X_\lambda^5$ and $Y_\lambda^5$ and of $X_1$ and $X_2$ are mirror.
\end{remark}
\subsection{Some Quotients by Subgroups of $\mathfrak{A}_{5}\rtimes H_4$}

\subsubsection{The group $\mathfrak{A}_{5}\rtimes H_4$}\label{quozientone}

Let $\varphi: {\mathfrak A}_5 \rightarrow Aut(H_5)$ be the homomorphism which maps an element $\tau \in {\mathfrak A}_5$ to the automorphism $\varphi(\tau)$ of $H_5$ such that $\varphi(\tau)(h)=\tau^{-1}h\tau$ in $Aut(X^5_{\lambda})$. An easy computation shows that $\varphi(\tau)(h)$ is the diagonal matrix in $H_5$ which is obtained from $h$ by applying the permutation $\tau$ to the diagonal entries of $h$. This allows one to compute the conjugacy classes of $\mathfrak{A}_{5}\rtimes H_4$. There are $25$ conjugacy classes. A representative has the form $(\tau, D)$, where $\tau \in {\mathfrak A}_5$ and $D$ is a diagonal matrix in $H_5$.  Below, we list the representatives with non-trivial fixed locus:
$$
id, ((12)(34), I), ((123), I),  $$
$$
(id, diag(1, 1,1, \alpha, \alpha^4)), \quad (id, diag(1, 1,1, \alpha^2, \alpha^3)),
$$
$$
(id, diag(1, 1,\alpha, \alpha, \alpha^3)), \quad (id, diag(1, 1,\alpha^2, \alpha^2, \alpha)),  $$
$$ ((12)(34), diag(1, 1,1,\alpha, \alpha^4)), \quad ((12)(34), diag(1, 1,1,\alpha^2, \alpha^3)),
$$
$$
((123), diag(1, 1,\alpha, \alpha, \alpha^3)), \quad ((123), diag(1, 1,\alpha^2, \alpha^2, \alpha)), $$
$$
((123), diag(1, \alpha,\alpha^2, \alpha^3, \alpha^4)), \quad (id, diag(1, \alpha^2,\alpha, \alpha^3, \alpha^4)).
$$

For each representative we must compute the fixed locus.  For instance, the fixed locus of  $((12)(34), I)$ is a curve of genus $6$ and a line contained in $X^5_{\lambda}$.  Near the fixed locus, we need to know the age of all representatives. This follows directly from \eqref{agegginverse} and the codimension of the fixed locus. For $((12)(34), I)$ the age is $1$. If we take into account the element $((123), I)$, the fixed locus has three components, a curve of genus $6$ and two points. The age of $((123), I)$ near the curve is $1$; near a point is $1$ and near the other fixed point is $2$.

The centralizer of an element can be computed from the size of the conjugacy class. For instance, the centralizer of $s=((12)(34), I)$ has order $20$ and fits into the short exact sequence:
$$
1 \rightarrow \Z/5\Z \rightarrow C_s \rightarrow \Z/2\Z \oplus  \Z/2\Z \rightarrow 1.
$$

By direct inspection, the quotient of the fixed locus of $ ((12)(34), I)$ by $C_s$ is the union of a rational curve and a genus two curve.

If we sum up all contributions according to the formula given by the orbifold cohomology, we finally get
$$
h^{1,1}=15, \quad h^{1,2}=5.
$$

This matches with previous results in \cite{St}, where the quotient is taken in two steps. First, note that $H_4$ is normal in $\mathfrak{A}_{5}\rtimes H_4$. The quotient of $X^5_{\lambda}$ is the singular variety $Y^5_{\lambda}$, which is the mirror as a toric variety. Next, lift the action of ${\mathfrak A}_5$ to a crepant resolution of $Y^5_{\lambda}$, and take the quotient by the lifting of all the automorphisms of ${\mathfrak A}_5$. This yields a singular variety that has a crepant resolution, which is actually a crepant resolution of $X^5_{\lambda}/\mathfrak{A}_{5}\rtimes H_4$.

\subsubsection{The Group $G=\Z/10\Z$}\label{sec: Z10} Let us consider the group $ \Z/10\Z\simeq G=\langle h_{(0,0,1,1,3)}, (12)(34)\rangle$. We observe that $G$ is generated by the element $$g_{10}:(x_1:x_2:x_3:x_4:x_5)\mapsto (x_2:x_1:\xi_5 x_4:\xi_5 x_3:\xi_5^3x_5).$$
The fixed locus of  $g_{10}$ consists the 2 points $(1:-1:0:0:0)$, $(0:0:1:-1:0)$, the fixed locus of $g_{10}^2$ consists of 10 points and the fixed locus of $g_{10}^5$ consists of the plane curve $C:=V(2x_1^5+2x_3^5+x_5^5-5\lambda x_1^2x_3^2x_5)$ and  of the line $(1:-1:s:-s:0)$. As in proof of Proposition \ref{cor: Hodge numbers quotient} we do not really need to know the linearization of $g_{10}$ near its fixed components, but only to observe that $age(g_{10}^i)+age(g_{10}^{-i})=3$ for $i=1,\ldots 9$, $i\neq 5$ and $age(g_{10}^5)=1$.\\
In order to compute the orbifold cohomology we observe that $g_{10}$ generates the centralizer both of $g_{10}^2$ and of $g_{10}^5$. A direct computation shows that $C/g_{10}$ is a smooth curve of genus 2, $Fix_{g_{10}^2}/g_{10}$ consists of 2 points.\\
We recall that $g_{10}$ acts trivially on the generator of the Picard group, and thus on $H^2(X_\lambda^5)$.
By \eqref{formula: rk trivial on cohom}, $p_{1,2}=9$ and thus, for any crepant resolution, $\widetilde{X_\lambda^5/g_{10}}$, of the quotient $X_\lambda^5/g_{10}$ we find:
$$h^{1,1}(\widetilde{X_\lambda^5/g_{10}})=1+4+4+2=11,
h^{2,1}(\widetilde{X_\lambda^5/g_{10}})=9+2=11.$$

\subsubsection{The Group $G=\Z/15\Z$}\label{sec: Z15} Let us consider the group $\Z/15\Z\simeq G=\langle h_{(0,0,0,1,4)}, (123)\rangle$. We observe that $G$ is generated by the element $g_{15}:(x_1:x_2:x_3:x_4:x_5)\mapsto (x_2:x_3:x_1:\xi_5 x_4:\xi_5^4x_5)$\\
The fixed locus of  $g_{15}$ consists the 2 points $P_1:=(1:\xi_3:\xi_3^2:0:0)$, $P_2:=(1:\xi_3^2:\xi_3:0:0)$, the fixed locus of $g_{15}^3$ consists of the curve of genus 6, $C:=V(x_1^5+x_2^5+x_3^5)$ and $P_i\in C$,  $g_{15}^5$ consists of the curve of genus 6, $D:=V(3x_1^5+x_4^5+x_5^5-5\lambda x_1^3x_4x_5)$ and of the two isolated points $P_i$. In order to compute the orbifold cohomology we observe that both $C/g_{15}$ and $D/g_{15}$ are smooth curve of genus 2.\\
Thus,  for any crepant resolution, $\widetilde{X_\lambda^5/g_{15}}$, of the quotient $X_\lambda^5/g_{15}$ we find
$$h^{1,1}(\widetilde{X_\lambda^5/g_{15}})=1+16=17,
h^{2,1}(\widetilde{X_\lambda^5/g_{10}})=9+12=21.$$

\subsubsection{The Group $G=\mathcal{D}_5$} Let us consider the dihedral group $\mathcal{D}_5\simeq G=\langle h_{(1,4,0,0,0)}, (12)(34)\rangle$ of order $10$. There are $4$ conjugacy classes in $\mathcal{D}_5$, namely:
$$
\{id\}, \{h_{(1,4,0,0,0)}, h_{(1,4,0,0,0)}^4\}, \{h^2_{(1,4,0,0,0)}, h_{(1,4,0,0,0)}^3\}, \{(12)(34)h^i_{(1,4,0,0,0)}\}_{i=1,\ldots 5}.$$
 The fixed locus of $h^i_{1,4,0,0,0}$, $i=1,2,3,4$ consists of the curve of genus 6, $C:=V(x_3^5+x_4^5+x_5^5)$ and the fixed locus of $(12)(34) h_{(1,4,0,0,0)}^j$ , $j=0,1,2,3,4$ consists of the curve of genus 6, $D_j:=V(2x_1^5+2x_3^5+x_5^5+5\lambda\xi^j x_1^2x_3^2x_5)$ and the line $l_j:=(t:-\xi^jt:s:-s:0)$. The centralizer group of every representative $s$ of a conjugacy class in $\mathcal{D}_5$ is generated by $s$ and thus its action on $Fix_s(X_\lambda^5)$ is trivial. We obtain $$h^{1,1}(\widetilde{X_\lambda^5/G})=1+4=5,  h^{2,1}(\widetilde{X_\lambda^5/G})=15+18=33.$$

\subsubsection{The Group $G=(\Z/5\Z)^2$} Let us consider the group $(\Z/5\Z)^2\simeq G=\langle h_{(0,1,2,3,4)}, (12345)\rangle$. Since both $h_{(0,1,2,3,4)}$ and $(12345)$ are fixed points free the quotient $X_\lambda^5/G$ is smooth and its Hodge numners are $h^{1,1}(X_\lambda^5/G)=1$, $h^{2,1}(X_\lambda^5/G)=5$ and the fundamental group of $X_\lambda^5$ is $(\Z/5\Z)^2$.\\

{\bf Acknowledgments.} The authors would like to thank Bert van Geemen for helpful suggestions. This work was partially supported by MIUR and GNSAGA.

\end{document}